\pdfoutput=1
\documentclass[11pt]{article}
\usepackage{amssymb}
\usepackage{graphicx}
\usepackage{tikzsymbols} 
\usepackage{amsmath}
\newtheorem{lm}{Lemma}

\newtheorem{prop}{Proposition}

\newtheorem{corol}{Corollary}

\newtheorem{ex}{Example}

\newtheorem{re}{Remark}

\newcommand{\ol}{\overline}

\newcommand{\ra}{\rightarrow}
\newcommand{\Ra}{\Rightarrow}

\newcommand{\s}{\subseteq}
\newcommand{\es}{\emptyset}

\newcommand{\LRa}{\Leftrightarrow}

\date{}

\begin{document}

\title{Compression with wildcards: Enumerating specific induced subgraphs, and packing them as well}

\author{Marcel Wild, Dept. of Mathematics, Stellenbosch University, South Africa}

\maketitle

\begin{abstract}
    Various algorithms have been proposed to enumerate all connected induced subgraphs of a graph $G=(V,E)$. As a variation we enumerate all "packings of connected sets", i.e. partitions $\Pi$ of $V$ with the property that each part of $\Pi$ induces a connected subgraph. More generally, for various types $T$ of graphs we do (one or both of) the following: (i) enumerate all type $T$ (induced) subgraphs of a given graph $G$, or (ii) enumerate all packings of type $T$ subgraphs of $G$.
\end{abstract}

\vspace{5mm}

\section{Introduction}

In recent articles the author has focused on  {\it numerically evaluating}  certain enumeration algorithms, while neglecting their theoretic assessment. The present article is the opposite: Eight Theorems state the {\it polynomial total time} performance of relevant algorithms. Most of these are not (yet) numerically evaluated. 
Always $G=(V,E)$ will be a simple graph with vertex set $V$ and $E$ is its set of edges $\{u,v\}$. (If $u,v\in V$ are such that $\{u,v\}\in E$, then $u,v$ are {\it adjacent} vertices.) 
Here comes the key concept.
For each (nonempty) $X\s V$  the subgraph $G[X]$ {\it induced} by $X$ has vertex-set $X$ and "as many edges it can get". Formally $G[X]:=(X,E(X))$ where
$$ E(X):=\Bigl\{\{u,v\}\in E:\ \{u,v\}\s X\Bigr\}.$$

\noindent
Our article  splits in two halves. The first  half is about enumerating various combinatorial objects individually, the second half "packs" these objects (in a sense defined in a minute).

\vspace{2mm}
As to the first half,   for eight types $T$ of graphs we do the following. Given $G=(V,E)$, we enumerate all $X\s V$ such that $G[X]$ is of type $T$. Moreover the enumeration takes place in a compressed format. These are the eight types (or properties):
$${\it connected,\  metric,\  convex,\  cycle\!\!-\!\!free,\  triangle\!\!-\!\!free,\  chordal,\  bipartite,\  chordless}$$
\noindent
The first three properties are dealt with in Section 5; here "metric" and "convex"  (two subtypes) are  certain exquisite ways of being connected.

The next four properties $T$ (Section 6) are such that $G[X]$ has property $T$ iff $X$ {\bf does not contain} chordless cycles of a specific type. Thus $G[X]$ being cycle-free (= being a forest) amounts to $X$ not containing {\it any} chordless cycles. And $G[X]$ being chordal, respectively bipartite, amounts to forbid chordless cycles of {\it length $\ge 4$}, respectively of {\it odd length}. 

The last property (Section 7) is an antipode to the four middle forur properties in that $G$ (often) {\bf does contain} chordless cycles. Specifically we call $G$ {\it chordless} if its connected components are either chordless paths or chordless cycles.

\vspace{2mm}

As to the second half of our article, given $G=(V,E)$ and a partition $\Pi=\{V_1,...,V_t\}$ of $V$ we say that $\Pi$ is a {\it packing of type $T$ graphs} if each graph $G[V_i]$ is of type $T$. In this scenario only two (not eight) types $T$ are investigated, but each with more effort. In Section 8 we generate all Cli-Pacs of $G$ (i.e. all packings of cliques), and in Section 9 all Conn-Pacs (i.e. all packings of connected sets).

\vspace{2mm}

The Sections 2,3,4,10 are of a more technical kind and will be surveyed in more detail at the beginning of the respective Sections. This much for now.  Given a powerset ${\cal P}(W)$ of some set $W$ and a set system ${\cal S}\s {\cal P}(W)$ we say that $X\in {\cal P}(W)$ is a {\it noncover} wrt $\cal S$ if $X\not\supseteq A$ for all $A\in {\cal S}$. Section 2 surveys the noncover $n$-algorithm of [W1] which produces all noncovers  wrt $\cal S$. This algorithm (and some variant, the implication $n$-algorithm) is the core ingredient to most of our results. Section 3 relates the matter to Boolean functions (specifically: Horn functions), and Section 10 introduces novel variants of the noncover $n$-algorithm. Finally Section 4 reviews basic facts about chordless paths, chordless cycles, and geodesics of a graph.

We will use the acronyms iff (= if and only if), wrt (=with respect to) and wlog (=without loss of generality).

\section{The noncover-algorithm and its variants: Part 1}

In a nutshell, the {\it noncover n-algorithm} of [W1] produces the set\footnote{Here "Mod" is an acronym for the modelset of a Boolean function; more details follow in Section 3. The $n$ in $n$-algorithm refers to the $n$-wildcard introduced below.}  $Mod({\cal S})\subseteq {\cal P}(W)$ of all noncovers  wrt ${\cal S} \subseteq {\cal P}(W)$ as follows.
Suppose ${\cal S}=\{A_1,A_2,...,A_h\}$. Starting with the powerset $\mbox{Mod}_0 = {\cal P}(W)$ we put
$$\mbox{Mod}_i : = \{X \in {\cal P}(W) : \ A_1, \cdots, A_i \not\subseteq X\}.$$
Provided $\mbox{Mod}_i$ is suitably encoded, one can smoothly sieve $\mbox{Mod}_{i+1}$ from $\mbox{Mod}_i$ while preserving the same kind of encoding. In the end $Mod_h=Mod({\cal S})$ is the family of all noncovers. Always $Mod({\cal S})$ is a set-deal of ${\cal P}(W)$  in the usual sense that $(\forall Y\in {\cal S})(\forall X\in {\cal P}(W))\ (X\s Y\Ra X\in {\cal S})$.

\vspace{3mm}
{\bf 2.1} To fix ideas let $W=\langle  9\rangle:=\{1,2,...,9\}$ and ${\cal S}:=\{A_1,A_2,A_3,A_4\}$, where
$$(1)\quad A_1:=\{1, 2, 4, 5\},\ A_2:=\{1, 2, 4, 7, 8, 9\},\ A_3:=\{2,5,8,9,\},\ A_4:=\{2,3,6, 9\}.$$
In the sequel we identify subsets of $W$ with their characteristic vectors $X\in\{0,1\}^9$ in the usual way and introduce two gadgets.
 One is the standard\footnote{In the literature often "$\ast$" is used instead of "2".} don't care symbol "2" which can freely be chosen as $0$ or $1$. The other is the wildcard $(n,n, \cdots, n)$ which means ``at least one $0$ here''. In other words, only $(1 1, \cdots, 1)$ is forbidden. Thus

$$\mbox{Mod}_1 = \{X \in {\cal P}(W): \ X \not\supseteq A_1 \} = (n,n,2,n,n,2,2,2,2),$$

\noindent 
which is row $r_1$ in Table 1. The acronym $PC = 2$ means that the {\it pending constraint}  to be imposed is the 2nd one. Thus we need to represent the family ${\cal F}\s r_1$ of all $X\in r_1$ that satisfy $X\not\supseteq A_2=\{1,2,4,7,8,9\}$. With respect to $A_1\cap A_2=\{1,2,4\}$ there are two types of such $X$'s:
\begin{itemize}
    \item[(i)] either $X\not\supseteq\{1,2,4\}$,  or (ii) $X\supseteq\{1,2,4\}$.
\end{itemize}

\noindent
A moment's thought shows that $r_3$ contains exactly the $X$'s of type (i), and $r_2$ the $X$'s of type (ii). (As to terminology, $r_2$ is a typical {\it 012n-row};  generally not all of  $0,1,2,n$ need to show up.)

\vspace{5mm}
\begin{tabular}{l|c|c|c|c|c|c|c|c|c|l} 
& 1 & 2 & 3 & 4 & 5 & 6 & 7 & 8 & 9 & \\ \hline 
&& & & & & & & & & \\
\hline 
$r_1=$ & $n$ & $n$ & 2 & $n$ & $n$ & 2 & 2 & 2 & 2 & $PC =2$ \\ \hline 
& & & & & & & & & & \\
\hline 
$r_3=$ & ${\bf n}$ & ${\bf n}$ & 2 & ${\bf n}$ & 2 & 2 & 2 & 2 & 2 & $PC = 3$ \\ \hline 
$r_2=$ & ${\bf 1}$ & ${\bf 1}$  & 2 & ${\bf 1}$ & 0 & 2 & $n$ & $n$ & $n$ & $PC = 4$\\ \hline
& & & & & & & & & & \\ \hline
$r_5=$ & 2 & ${\bf 0}$ & 2 & 2 & 2 & 2 & 2 & 2 & 2 & final\\ \hline
$r_4=$ & $n_1$ & ${\bf 1}$ & 2 & $n_1$ & $n_2$ & 2 & 2 & $n_2$ & $n_2$ & $PC = 4$\\ \hline
$r_2=$ & 1 & 1 & 2 & 1 & 0 & 2 & $n$ & $n$ & $n$ & $PC = 4$\\ \hline
& & & & & & & & & & \\
 \hline
$r_7=$ & $n_1$ & 1 & 2 & $n_1$ & 2 & 2 & 2 & 2& ${\bf 0}$ & final \\ \hline
$r_6=$ & $n_1$ &  1& $n_3$ & $n_1$ & $n_2$ & $n_3$ & 2 & $n_2$ & ${\bf 1}$ & final\\ \hline
$r_2=$ & 1 & 1 & 2 & 1 & 0 & 2 & $n$ & $n$ & $n$ & $PC = 4$ \\ \hline
& & & & & & & & & & \\ \hline
$r_9=$ & 1 & 1 & 2 & 1 & 0 & 2 & 2 & 2 & ${\bf 0}$ & final \\ \hline
$r_8=$ & 1 & 1 & $n_2$ & 1 & 0 & $n_2$ & $n_1$ & $n_1$ & ${\bf 1}$ & final \\ \hline
\end{tabular}

\vspace{3mm}
{\sl Table 1: The noncover $n$-algorithm; snapshots of its working stack}

\vspace{4mm}
  By construction $r_2$ and $r_3$ satisfy the second constraint, and incidentally $r_2$ also satisfies the third (i.e. $A_3 \not\subseteq X$ for all $X \in r_2$). Hence $r_2$ has $PC = 4$, and $r_3$ has $PC=3$.
The {\it working stack} now is $\{r_2, r_3\}$ with $r_3$ being on top. We keep on picking the top row $r$, impose its pending constraint, and thereby replace $r$ by its sons. If a top row has no more pending constraints, it is {\it final} and gets moved to a save place. Proceeding in this manner the working stack $\{r_2, r_3\}$ becomes $\{r_2, r_4, r_5\}$. Upon removing its final top row $r_5$ we get $\{r_2, r_4\}$ which leads to $\{r_2, r_6, r_7\}$. Upon removing $r_6, r_7$ we get $\{r_2\}$ which yields the final rows $r_8, r_9$.
Since the final rows are mutually disjoint, the set $Mod_6=Mod({\cal S})$ of noncovers wrt $\cal S$ has cardinality
$$(2)\quad |Mod({\cal S})|=|r_5| + |r_6| + |r_7| + |r_8| + |r_9| = 2^8 + 2 \cdot 3^3 + 2^5 \cdot 3 + 3^2 + 2^4 = 431.$$
Theorem $1'$ below is a special case of Theorem 1 in Section 3.

\vspace{3mm}

{\bf Theorem $1'$: }{\it The noncover $n$-algorithm calculates the $N$ noncovers wrt $\{A_1,...,A_h\}\s W$ in time $O(Nh^2|W|^2)$.}

\vspace{3mm}

{\bf 2.2}  Let $W$ be any set and $(A,B)\in {\cal P}(W)\times{\cal P}(W)$ an ordered pair of subsets. We henceforth call $(A,B)$ an {\it implication}, but prefer to write it as $A\to B$. One calls $A$ is  the {\it premise}  and $B$ the {\it conclusion} of the implication.
If $\Sigma=\{A_1\to B_1,...,A_h\to B_h\}$ is any  implication-family {\it on W} (i.e. all $A_i,B_i\s W$), one says that $Y\s W$ is $\Sigma${\it -closed} if for all $1\le i\le h$ it holds that $(A_i\not\s Y\ or\ B_i\s Y)$. Put another way, $A_i\s Y$ {\it implies} $B_i\s Y$.

Our first variant of  the noncover $n$-algorithm (others follow in Section 9) is the most widely applicable one. It  is the {\it implication $n$-algorithm} and renders
 the family $Mod(\Sigma)$ of all $\Sigma$-closed sets. Similar to the sets $A_i\in{\cal S}$ in 2.1, now the implications $(A_{i}\to B_{i})\in\Sigma$ get imposed one-by-one. Thus suppose the row $r$ below satisfies the first $i$ implications and now the $(i+1)$-th implication
$A_{i+1}\to B_{i+1}$, say $\{2,3\}\to \{7,9\}$, needs to be imposed:

\vspace{2mm}
\begin{tabular}{l||c|c||c|c|c|c|c|c|c|l} 
& 1 & 2 & 3 & 4 & 5 & 6 & 7 & 8 & 9 & \\ \hline 
&& & & & & & & & & \\
\hline 
$r=$ & $n_1$ & $n_1$ & $n_2$ & $n_2$ & $n_2$ & 0 & 2 & 2 & 2 & $PC\quad  is\quad i+1$ \\ \hline 
&& & & & & & & & & \\ \hline
$\rho_1=$ & 2 &  {\bf 0} & $n_2$ & $n_2$ & $n_2$ & 0 & 2 & 2 & 2 &  \\ \hline
$\rho_2=$ & 0 &  {\bf 1} &  {\bf 0} & 2 & 2 & 0 & 2 & 2 & 2 &  \\ \hline
$\rho_3=$ & 0 &  {\bf 1} &  {\bf 1} & $n_2$ & $n_2$ & 0 & {\bf 1} & 2 &  {\bf 1} &  \\ \hline
\end{tabular}
\vspace{2mm}

{\it Table 2: Glimpsing the implication $n$-algorithm}

\vspace{2mm}
One checks that $\rho_1\uplus \rho_2$ is the set of $X\in r$ that {\it trivially satisfy} $A_{i+1}\to B_{i+1}$ in the sense that $A_{i+1}\not\s X$ (so $X$ is a noncover of $A_{i+1}$). In contrast $\rho_3$ contains those $X\in r$ that
{\it nontrivially satisfy} $A_{i+1}\to B_{i+1}$. This means that $A_{i+1}\s X$, and thus by definition  $B_{i+1}\s X$. Hence  $\rho_1\uplus \rho_2\uplus \rho_3$ contains  those $X\in r$ that satisfy $A_{i+1}\to B_{i+1}$. The next action in the overall implication $n$-algorithm would cancel the top row $r$ of the stack  (cf Table 1) and replace it by the rows
$\rho_1,\rho_2, \rho_3$. Each one of these rows has a pending constraint (PC) $\le i+2$.

As is well known, if $\Sigma$ is any implication-family, then the set system
$Cl(\Sigma)$ of all $\Sigma$-closed sets is a {\it closure system}, i.e. $W\in Cl(\Sigma)$ and from $X,Y\in Cl(\Sigma)$ follows $X\cap Y\in C(\Sigma)$.
Like Theorem $1'$ also Theorem $1''$ is a consequence of Theorem 1 in Section 3.

\vspace{2mm}
{\bf Theorem $1''$: }{\it If $\Sigma:=\{A_{1}\to B_{1},..,A_{h}\to B_{h}\} $ is an implication-family on the set $W$, then the implication $n$-algorithm produces the $N$ many $\Sigma$-closed sets in time $O(Nh^2|W|^2)$.}

\section{About  Horn clauses and Horn CNFs}

 It is time to relate all of the above to Boolean functions. Thus recall  [CH] that a {\it clause}
is  a disjunction of literals, such as $C:=\ x_2\vee \ol{x}_4 \vee \ol{x}_6\vee x_7$. The first and last literal are {\it positive}, the middle two are {\it negative}. A conjunction $f:=C_1\wedge C_2\wedge\ldots\wedge C_k$ of clauses $C_i$ is a {\it Conjunctive Normal Form (CNF)}.

By definition a {\it Horn   clause} has at most one positive literal, so $C$ above is {\it not} a Horn   clause.
A conjunction of Horn clauses is called {\it Horn CNF}. Not all Horn CNFs are satisfiable, an obvious example being $x_1\wedge \ol{x}_1$. As opposed to arbitrary CNFs the satisfiability of a Horn CNF $f$ can be tested in polynomial time (details in 3.1).

Let $f$ be a {\it Boolean function}, i.e. a function of type $f:\{0,1\}^n\to\{0,1\}$  (such $f$ can e.g. be defined by a CNF in obvious ways). Any bitstring $y\in\{0,1\}^n$ with $f(y)=1$ is a {\it model} of $f$; the set $Mod(f)$ of all models is the {\it modelset of $f$}. For instance consider this Horn CNF:
$$(3)\quad f:=(\ol{x}_1\vee \ol{x}_2\vee\ol{x}_4\vee\ol{x}_5)\wedge (\ol{x}_1\vee\ol{x}_2\vee\ol{x}_4\vee\ol{x}_7\vee\ol{x}_8\vee\ol{x}_9)\wedge 
(\ol{x}_2\vee \ol{x}_5\vee\ol{x}_8\vee\ol{x}_9)\wedge (\ol{x}_2\vee \ol{x}_3\vee\ol{x}_6\vee\ol{x}_9).$$
It is clear that $Mod(f)$ matches $Mod({\cal S})$ from  2.1. For instance $(1,1,1,0,0,0,1,1,1)\in Mod(f)$ matches
$\{1,2,3,7,8,9\}\in Mod({\cal S})$. Generally the noncover $n$-algorithm yields exclusively  {\it negative}  clauses (i.e. having only negative literals).

\vspace{3mm}

Consider now implications, such as $\{2,3\}\ra\{7,9\}$. It amounts to ($\{2,3\}\ra\{7\}$ and $\{2,3\}\ra\{9\}$), and this  translates to the Horn CNF $(\ol{x}_2\vee \ol{x}_3\vee x_7)\wedge (\ol{x}_2\vee \ol{x}_3\vee x_9)$. It follows that  for each implication-family $\Sigma$ there is a Horn CNF $f$ such that $Mod(\Sigma)=Mod(f)$.
Generally the implication $n$-algorithm yields exclusively  {\it pure} Horn clauses (i.e. having a positive literal). 

\vspace{3mm}
{\bf 3.1} The {\it Horn $n$-algorithm} accepts arbitrary Horn CNFs and according to [W1, Thm.2] performs as follows.

\vspace{3mm}
{\bf Theorem 1: }{\it Let $f$ be a satisfiable Horn CNF with $h$ clauses and altogether $w$ literals. Then the $N$ models in $Mod(f)\s \{0,1\}^w$  can be enumerated in time $O(Nh^2w^2)$.}

\vspace{2mm}
The Theorems $1',1'',1$ will be invoked in the proofs of Theorems 3,4,5,6,7 and 2,11 and 9 respectively. Albeit Theorem 1 would suffice in all instances, it is worthwile to point out the various levels of sophistication. Not least because the noncover $n$-algorithm and the implication $n$-algorithm are more susceptible to taylor-made adaptions than the Horn $n$-algorithm (as seen in Section 10).

\vspace{2mm}
If $f$ in Theorem 1 was an unrestricted Horn CNF then the bound $O(Nh^2w^2)$ would be false by the following reason. If $N=0$ then $O(Nh^2w^2)=0$, yet it takes time $>0$ (in fact time $O(hw)$) to detect that $Mod(f)=\es$, i.e. that $f$ is insatisfiable. But then again, the statement of Theorem 1 (as it is) suffices\footnote{Notice that $N>0$ is automatically satisfied in Theorem $1'$ and Theorem $1''$. Indeed, in the first case $\es\in Mod({\cal  S})$, in the second case
$W=(1,1,...,1,1)\in Mod(\Sigma)$.} since in our upcoming applications $f$ is always satisfiable.

\vspace{2mm}
{\bf 3.2} Let us take the satisfiability issue further in another direction. Namely, in some applications it may be desirable not to compute  the whole of $Mod(f)$, but only the subfamily
$$Mod(f,\ge k):=\{y\in Mod(f):\ |ones(y)|\ge k\}$$
for some suitable $k>0$. Although $Mod(f,\ge k)$ cannot be enumerated in polynomial total time one can  cut (but not prevent) the useless production of duds, i.e  of intermediate 012n-rows $r$ which are such that no successor row of $r$ is in $Mod(f,\ge k)$.
To fix ideas consider the Horn CNF
$$F:=\ (\ol{x}_1\vee x_2)\wedge (\ol{x}_2\vee\ol{x}_3\vee\ol{x}_6)\wedge (x_3\vee\ol{x}_4\vee\ol{x}_7)\wedge (x_1\vee \ol{x}_5).$$
For the sake of clarity we stick to 012-rows in this example. Thus suppose the Horn 012-algorithm found that $r:=(2,2,2,2,0,1,1)$ is feasible, i.e. $r\cap Mod(F)\neq\es$. What if this is not enough and we accept $r$ (and process it further) only if $r\cap Mod(F,\ge 4)\neq\es$? First observe that $r\cap Mod(F)=Mod(F_r)$, where $F_r$ is obtained from $F$ by setting $x_5:=0,\ x_6:=1,\ x_7:=1$. Thus
$$F_r=\  (\ol{x}_1\vee x_2)\wedge (\ol{x}_2\vee\ol{x}_3)\wedge (x_3\vee\ol{x}_4).$$ 
If there are $y\in Mod(F_r)$ with (say) $y_1=y_4=1$ then $y\in r\cap Mod(F,\ge 4)$, and so we could accept $r$. To decide the existence of such $y$'s put $x_1=x_4:=1$ in $F_r$. This yields 
$$F_r^{1,4}:=\ x_2\wedge (\ol{x}_2\vee\ol{x}_3)\wedge x_3,$$
which is satisfiable. On the other hand, e.g. $F_r^{1,3}:=\ x_2\wedge \ol{x}_2\wedge 1$ is insatisfiable. One checks that exactly 3 out of $\binom{4}{2}=6$ pairs $\{i,j\}\s twos(r)=\{1,2,3,4\}$ are "good" and make $F_r^{i,j}$ satisfiable. Therefore we accept $r$ and "process" it by substituting it with its three good "sons" $r^{i,j}$ (which, to spell it out, are defined by $zeros(r^{i,j}):=zeros(r),\ ones(r^{i,j}):=ones(r)\cup\{i,j\},\ twos(r^{i,j}):=twos(r)\setminus\{i,j\}$).

\vspace{2mm}
Generally one may wish to "look $t$ steps ahead"  in the sense that during the algorithm all feasible rows $r$ are processed as follows. Among the $\binom{|twos(r)|}{t}$ "candidate sons" of $r$ pick those (if any) that satisfy
$|ones(r')|=|ones(r)|+t$ and substitute $r$ with them; this entails the case that $r$ gets cancelled without substitutes.
As is plausible from above, one can show that trimming the Horn $n$-algorithm in this manner pushes the cost from $O(Nh^2w^2)$ to $O\big(N(h^2w^2+\binom{w}{t}hw)\big)=O(Nh^2w^t)$.

\vspace{2mm}

{\bf 3.3} Recall that Horn CNFs with exclusively negative clauses give rise to the $n$-wildcard (whose meaning was "at least one 0 here"). If these negative clauses display specific  patterns then the use of additional wildcards is beneficial. We start with  informal definitions:

\begin{itemize}
    \item[(4)] $(n(2),n(s),...,n(2))$ means "at least two 0's here".
     \item[] $(\epsilon(2),\epsilon(2),...,\epsilon(2))$ means "at most two 1's here".
     \item[] $(a,\epsilon(2),\epsilon(2),...,\epsilon(2))$ means "if $a=1$ then at most two 1's on the rest".
      \item[] $(\epsilon,\epsilon,...,\epsilon)$ means "at most one 1 here".
       \item[] $(a,\epsilon,\epsilon,...,\epsilon)$ means "if $a=1$ then at most one 1 on the rest".
     \item[] $(a,c,c,...,c)$ means "if $a=1$ then only 0's on the rest".
      \item[] $(\gamma,\gamma,...,\gamma)$ means "exactly one 0 here".
     \end{itemize}

   Hence each wildcard in (4) describes a set of bitstrings by way of specific property that the bitstrings must possess. 
   For instance, $(\epsilon,\epsilon):=\{(0,0),(0,1),(1,0)\}$ and $(a,\epsilon,\epsilon):=(0,2,2)\cup (1,\epsilon,\epsilon)$.
   Here comes a  formal definition in terms of Horn CNFs. To avoid excessive notation we look at Horn CNFs on six literals.

     \vspace{3mm}

\begin{tabular}{|l|c|c|c|c|c|c|c|} 
1& 2 & 3 & 4 & 5 & 6 &    & defining Horn CNF   \\ \hline 
& & & & & &  &  \\ 
 $\epsilon(2)$ & $\epsilon(2)$ & $\epsilon(2)$ & $\epsilon(2)$ & $\epsilon(2)$ & 2 &    &
$(\ol{x_1}\vee \ol{x_2}\vee \ol{x_3})\wedge (\ol{x_1}\vee \ol{x_2}\vee \ol{x_4})
\wedge...\wedge (\ol{x_3}\vee\ol{x_4}\vee \ol{x_5})$    \\ \hline 
& & & & & &  &  \\  

 $\epsilon(2)$ & $\epsilon(2)$ & $\epsilon(2)$ & $\epsilon(2)$ & $\epsilon(2)$ & $a$ &    &
$\ol{x_6}\vee \big[(\ol{x_1}\vee \ol{x_2}\vee \ol{x_3})\wedge\ldots \wedge (\ol{x_3}\vee\ol{x_4}\vee \ol{x_5})\big]$    \\ \hline 
& & & & & &  &  \\ 
$\epsilon$ & $\epsilon$ & $\epsilon$ & $\epsilon$ & $\epsilon$ & 2 &    &
$(\ol{x_1}\vee \ol{x_2}\vee \ol{x_3}\vee\ol{x_4})\wedge\ldots\wedge (\ol{x_2}\vee \ol{x_3}\vee\ol{x_4}\vee \ol{x_5})$    \\ \hline 

& & & & & &  &  \\  
$\epsilon$ & $\epsilon$ & $\epsilon$ & $\epsilon$ & $\epsilon$ & $a$ &    &
$\ol{x_6}\vee\big[(\ol{x_1}\vee \ol{x_2}\vee \ol{x_3}\vee\ol{x_4})\wedge\ldots\wedge (\ol{x_2}\vee \ol{x_3}\vee\ol{x_4}\vee \ol{x_5})\big]$    \\ \hline 

& & & & & &  &  \\  
$c$ & $c$ & $c$ & $c$ & $c$ & $a$ &    &
$\ol{x_6}\vee\big[\ol{x_1}\wedge \ol{x_2}\wedge \ol{x_3}\wedge\ol{x_4}\wedge\ol{x_5}\big]$    \\ \hline 

\end{tabular}

\vspace{2mm}
{\sl Table 3: Formal definitions of some wildcards}

\vspace{2mm}
\noindent

A few remarks are in order. Let $f_1$ be the Horn CNF coupled to $(\epsilon(2),...,\epsilon(2),2)$ and consider any bitstring $y\in\{0,1\}^6$ with at least three 1's among its first five components. Then $y$ violates at least one of the $\binom{5}{3}$ Horn clauses. For instance $y':=(1,0,1,1,0,1)$ violates $\ol{x}_1\vee \ol{x}_3\vee \ol{x}_4$. Therefore $Mod(f_1)$ contains exactly the bitstrings  with at most two 1's among the first five components. 

Similarly one explains the Boolean formulas $f_2,...,f_5$ corresponding to the other wildcards in Table 3. As to $f_4$ and $f_5$, they are not themselves a Horn CNFs, but e.g. $f_5$ is equivalent to the Horn CNF $(\ol{x}_6\vee \ol{x}_1)\wedge\ldots\wedge (\ol{x}_6\vee \ol{x}_5)$. As to  $(n(2),...,n(2))$ and $(\gamma,...,\gamma)$ occuring in (4), they will be interwoven in Section 10.3 in a way that also matches a Horn CNF.

\vspace{2mm}
{\bf 3.4} In later Sections all wildcards in (4) (and of course the $n$-wildcard) will be seen in action in graph theoretic settings. Except for $(a,c,...,c)$, all of them are new. For a broader visions let us mention some previously used wildcards (which will not be used later on):

\begin{itemize}
    \item[(5)] $(e(k),e(k),...,e(k))$ means "at least $k$ many 1's here".
     \item[] $(e,e,...,e):=(e(1),e(1),...,e(1))$.
     \item[] $(g(k),g(k),...,g(k))$ means "exactly $k$ many 1's here".
      \item[] $(g,g,...,g):=(g(1),g(1),...,g(1))$.
     \item[] $(a,b,b,...,b)$ means "if $a=1$ then only 1's on the rest".
     \end{itemize}

     We chose the letter $\epsilon$ in (4) because it resembles $e$, and because "at most one 1" dualizes "at least one 1". Similarly $\gamma$ is the Greek version of $g$ and "exactly one 0 here" dualizes "exactly one 1 here".
     (It remains the question why\footnote{This is due to my German background: $e$ derives from Eins (=one), $n$ derives from Null (=zero), and $g$ derives from genau (=exactly).} the letters $e,n,g$ were chosen in the first place.)

\section{Chordless paths, chordless cycles, and geodesics}

To fix terminology, the graph $G_1=(V^1,E_1)$ in Figure 1 (never mind the dashed edges) has the $8-3$ path $R=(8,4,7,6,3)$ with underlying vertex set $V^1(R):=\{3,4,6,7,8\}$. Changing direction yields $R':=(3,6,7,4,8)$ and obviously $V^1(R')=V^1(R)$. The path $R$
 has the {\it chord} $\{8,7\}$, i.e. an edge that connects  non-consecutive vertices of the path. In contrast the path $P=(8,7,6,3)$ is a {\it chordless } path (clp). 
 Each clp (but not only them) is uniquely determined by its underlying vertex-set. Formally, for any graph $G=(V,E)$, if $X\s V$ is such that $X=V(P)=V(Q)$ for cl paths $P$ and $Q$, then $Q=P\ or\ Q=P'$. Put another way, a path $R$ in $G$ is chordless iff  the induced graph $G[V(R)]$ has as many edges as $R$.

 \begin{center}
\includegraphics[scale=0.72]{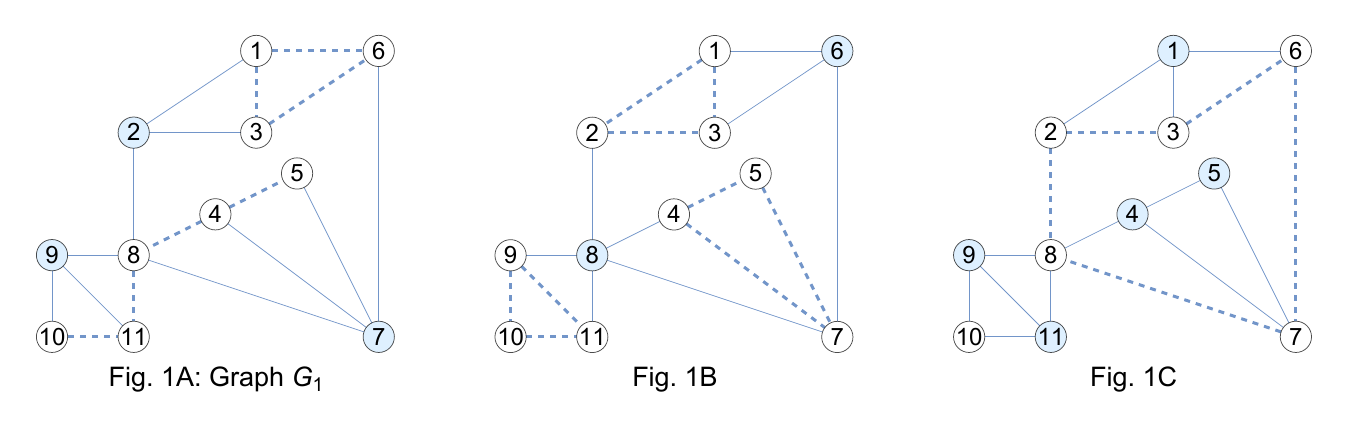}
\end{center}

 \vspace{3mm}

There is an even shorter path that leads from $8$ to $3$ in $G_1$, namely $(8,2,3)$. 
A path of shortest length among all paths from $s$ to $t\ (s\neq t)$ is called $s-t$ {\it geodesic}, or just {\it geodesic} if the endpoints are irrelevant.  Each geodesic is a cl path, but not conversely.  
We write $Geo(s,t)$ for the set of all $s-t$ geodesics, and $Geo(G)$ for the set of all geodesics of $G$. Similarly $CLP(s,t)$ and $CLP(G)$ are defined wrt chordless paths.

\vspace{2mm}
Since in our context cycles\footnote{Recall that cycles are circuits which do not repeat vertices and (hence) not edges. Observe that our terminology is different from the one in [S,p.20].} will never be oriented, we can and will formally identify each cycle with the set of edges that constitute it. As for paths,
different cycles  can  have the same underlying vertex-set; an extreme example is the complete graph  on $n$ vertices, where all of its $(n-1)!$ length $n$ cycles $C$ share the same vertex-set $V(C)=V$. Also like clps, chordless cycles (clcs) are uniquely determined by their vertex set.
Let $CLC(G)$ be the set of all chordless cycles of $G$. (Depending on the application we will identify a clc either with its underlying edge-set or its underlying vertex-set.)
 
In the remainder of Section 4 we sketch how to generate  $CLP(G), Geo(G)$, and $CLC(G)$.

\vspace{2mm}
{\bf 4.1} Consider again the graph $G_1$. One way to calculate $CLP(G)$ proceeds inductively as follows. Given that the set $CLP[k]$ of al cl paths of length $k$ has been calculated, scan all $P\in CLP[k]$  and extend $P$ on the right in all possible ways that yield a member $P'$ of $CLP[k+1]$. For instance $CLP[2]=\{ (1,2,8),(1,6,7),(2,1,6),...,(8,11,10)\}$. Then say $P=(2,1,6)$ cannot be extended to $(2,1,6,3)$ (because $\{1,3\}$ is a chord), but $Q=(2,1,6,7)$ is a new member of $CLP[3]$. As to "new member", recall that $(7,6,1,2)$ is considered the same as $(2,1,6,7)$. Therefore, listing only cl paths whose last vertex has higher value than the first, avoids repetitions and exhausts $CLP[k+1]$.

\vspace{2mm}
{\bf 4.2} It will often be suitable to  endow the  vertex-set  of $G=(V,E)$ with an arbitrary linear order $<$; consequently, if  $P$ is a $s-t$ geodesic and $s<t$, we may convene to traverse $P$ from $s$ to $t$.
The set $Geo(G)$ of all geodesics of $G$ can be calculated very much like $CLP(G)$ in 4.1. The necessary extra ingredient to determine the feasible extensions of geodesics $P\in Geo[k]$ is the $n\times n$ distance matrix $D(G)$. By definition its $(s,t)$-entry gives the common length of all $s-t$ geodesics, i.e. the {\it distance} between $s$ and $t$. Interestingly this algorithm, called {\tt Natural APAG} in [W3], often beats the standard\footnote{The standard procedure splits $Geo(G)$ into the parts $G(s,t)\ (s<t)$ and calculates all sets $Geo(s,t)$ with depth-first-search. This is more clumsy than the 4.1 way. Admittedly, if {\it only} $Geo(s,t)$ is required for two specific $s\neq t$, then the 4.1 way is not applicacable.} procedure.

\vspace{2mm}
{\bf 4.3} The calculation of the set $CLC(G)$ of all cl cycles of $G$ is more cumbersome. One naive option is to first calculate $CLP(G)$ (the  4.1 way or otherwise), and to investigate each $P\in CLP(G)$. If the start vertex of $P$ is adjacent to its end vertex, then $V(P)$ is the vertex set of a  cl cycle, otherwise not. Trouble is, most cl paths may not yield a cl cycle, and so the generation of $CLP(G)$ was overkill. Furthermore, if $P$ happens to yield (vertex-wise) a cl cycle $V(P)$ then $|V(P)|-1$ other members of $CLP(G)$ yield the {\it same} cl cycle; thus again a waste of time.
Several authors have proposed methods to circumvent these  issues. None of these attempts is as clear-cut as "the 4.1 way". One of the more reader-friendly ones is [DCLJ].

(In Section 7 we offer a clear-cut algorithm for a problem which is related insofar that chordless paths and chordless cycles simultaneously take the stage.)

Each  cycle $C\s V$ with $|C|=3$ is automatically  chordless and is called a {\it triangle}. Let $Triangles(G)\s CLC(G)$ be the family of all triangles of $G$. Other than $CLC(G)$, the calculation of $Triangles(G)$ is painless; more on that in Section 10.1.

\section{Induced subgraphs that are connected or metric or convex }

Here we  show how all connected, or all metric, or all convex  induced subgraphs of $G=(V,E)$ can be generated; in fact all of them  in compressed fashion.
By abuse of language we say that $X$ is {\it connected} when $G[X]$ is connected. Hence $X$ is connected 
iff for all $s\neq t$ in $X$ there is some $s-t$ path $P$ with $V(P)\s V$. 
Fortunately it suffices\footnote{If there is a chord in $P$, say $\{v_2,v_k\}\in E$, then $\{v_2,v_k\}\in E(X)$, and so $(s,v_1,v_2,v_k,t)$ is a path within $G[X]$ with at least one chord less than $P$. The claim follows by induction. } to demand the existence of a {\it chordless} path $P$ with $V(P)\s V$. This matters from an algorithmic point of view since there may be way less chordless $s-t$ paths
than ordinary $s-t$ paths. The definition in formulas:

\vspace{3mm}
\begin{itemize}
    \item[(6i)] $X\s V$ is {\it connected} iff  $\Big(\forall s,t\in X\Big)\Big( s\neq  t\Ra \exists P\in CLP(s,t)\ {\rm with}\ V(P)\s X\Big)$.
\end{itemize}

\noindent
Upon switching $\exists  P$ with $\forall P$, and switching  $CLP(s,t)$ with $Geo(s,t)$, one obtains altogether four related concepts. Specifically:

\vspace{3mm}
\begin{itemize}
    \item[(6ii)] $X\s V$ is {\it metric} iff  $\Big(\forall s,t\in X\Big)\Big( s\neq  t\Ra \exists P\in Geo(s,t)\ {\rm with}\ V(P)\s X\Big)$.
    \item[(6iii)] $X\s V$ is {\it ge-convex} iff  $\Big(\forall s,t\in X\Big)\Big( s\neq  t\Ra \forall P\in Geo(s,t)\ {\rm it\ holds\ that\ }V(P)\s X\Big)$.
    \item[(6iv)] $X\s V$ is {\it mo-convex} iff  $\Big(\forall s,t\in X\Big)\Big( s\neq  t\Ra \forall P\in CLP(s,t)\ {\rm it\ holds\ that\ }V(P)\s X\Big)$.
\end{itemize}

\noindent
Here "mo" is an acronym for "monophonically" and "ge" an acronym for "geodesically" (see [FJ],[W2]). 
It follows from  (6i) to (6iv) that these  implications take place:

$$(7)\quad mo\!\!-\!\!convex\ \Ra\  ge\!\!-\!\!convex\ \Ra\ metric\ \Ra\ connected$$

\noindent
These three sets of vertices in $G_1$ prove that none of the implications in (7) are reversible:
$\{2,3,8\},\ \{1,2,6\},\ \{3,6,7,8\}$. 

\vspace{2mm}
For given $G=(V,E)$ we let $MoConv(G)$ be the set of all mo-convex sets $X\s V$. Similarly
$GeConv(G), Metric(G)$, and $ Conn(G)$ are defined. The remainder of Section 5 is dedicated to a compressed enumeration of these four set families. They come in pairs; $MeConv(G),\ GeConv(G)$ in 5.1, and $Metric(G), Conn(G)$ in 5.2.
Numerical experiments follow in 5.3.

\vspace{3mm}
{\bf 5.1} For  $s\neq t$ in $V$ let $V(CLP(s,t))$ be the union of all sets $V(P)$ where $P$ ranges over $CLP(s,t)$. Evidently $V(CLP(s,t))=V(CLP(t,s))$ and  $\{s,t\}\in E\Ra V(CLP(s,t))=\{s,t\}$. With this in mind we define
$$W_G:=\{(s,t)\in V\times V:\ s<t\ and\ \{s,t\}\not\in E\}.$$ 
The definition $V(Geo(s,t))$ is analogous to the one of $V(CLP(s,t))$ and we put

$$(8)\quad \Sigma_{mo}^G:=\Big\{\{s,t\}\to V(CLP(s,t)):\ \{s,t\}\in W_G\Big\}\ and\ \
\Sigma_{ge}^G:=\Big\{\{s,t\}\to V(Geo(s,t)):\ \{s,t\}\in W_G\Big\}.$$

\vspace{2mm}
{\bf Theorem 2: }
\begin{itemize}
\item[(a)]
{\it If $\Sigma_{mo}^G$ in (8) is known, then $MoConv(G)$ can be calculated in time $O(N_1|V|^8)$, where $N_1:=|MoConv(G)|$.} 
\item[(b)]{\it If $\Sigma_{ge}^G$ in (8) is known, then $GeConv(G)$ can be calculated in time $O(N_2|V|^8)$, where $N_2:=|GeConv(G)|$.}
\end{itemize}

\vspace{2mm}
{\it Proof.} We only prove (a), the argument for (b) is similar. From $(6iv)$ and (8) follows that $MoConv(G)=Cl(\Sigma_{mo}^G)$. Putting $h:=|\Sigma_{mo}^G|$ it follows from Theorem $1''$ that $Cl(\Sigma_{mo}^G)$ can be calculated (in potentially compressed fashion) in time $O(N_1h^2|W|^2)=O(N_1\cdot(|V|^2)^2\cdot (|V|^2)^2)=O(N_1|V|^8)$. $\square$

\begin{center}
\includegraphics[scale=0.9]{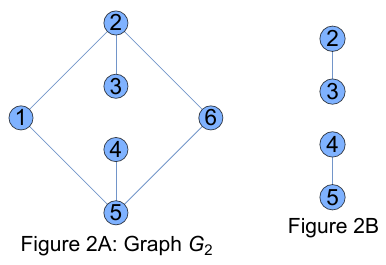}
\end{center}

\vspace{3mm}
{\bf 5.2} Concerning the pair $Conn(G)$ and $Metric(G)$, first observe that $Conn(G)$ is no set ideal of ${\cal P}(V)$, unless  $G$ is complete.  Usually $Conn(G)$ is not even  a closure system. For instance for $G_2$ in Fig.2A it holds  that  $X_1:=\{1,2,3,4,5\}$ and $X_2:=\{2,3,4,5,6\}$ belong to $Conn(G_2)$, yet this is not the case for $X_1\cap X_2$ (see Fig.2B).
It follows that generally (for exceptions see [W2]) neither the noncover $n$-algorithm nor the implication $n$-algorithm are fit to generate $Conn(G)$. In other words, whatever implication-family one may come up with, one gets $Cl(\Sigma)\neq Conn(G)$.

\vspace{2mm}
On a Boolean logic level,
what  makes the enumeration of  $Conn(G)$ more complicated than the one of $MoConv(G)$
is that $\forall P$  in (6iv)  becomes $\exists P$  in (6i). The way out is to view $Conn(G)$ as the modelset of some suitable Boolean function which is in CNF format and has very long clauses; see [W2]. 

\vspace{2mm}
Likewise, what  makes the enumeration of  $Metric(G)$ more complicated than the one of $GeConv(G)$
is that $\forall P$  in (6iii)  becomes $\exists P$  in (6ii). To summarize, what unites $Conn(G)$ and $Metric(G)$ is $\exists$, what separates them are the inputs $CLP(G)$ and $Geo(P)$.

\vspace{3mm}
{\bf 5.3} As to numerical experiments, the algorithm {\tt AllMetricSets} of [W2] enumerates $Metric(G)$.
Albeit polynomial total time cannot be proven, it is fairly efficient in practise.
Specifically, it was pit against the Mathematica command {\tt BooleanConvert} (option "ESOP") which  enumerates all models of {\it any} Boolean function (and, like the author, does so in a compressed fashion). {\tt BooleanConvert} is  disadvataged in that {\tt AllMetricSets} is taylored to the relevant type of Boolean function. On the other hand, {\tt BooleanConvert} is a built-in command (coded in C or C++), whereas {\tt AllMetricSets} is coded in high-level Wolfram Language\footnote{This is the only language I use since 30 years.}. This makes for an interesting competition. In fact,  {\tt AllGeConvexSets} (which calculates $GeConv(G)$) was  also coded in high-level Wolfram Language and participates in the competition. Recall from Theorem 2(b) that for  {\tt AllGeConvexSets} polynomial total time {\it can} be proven.

\vspace{3mm}
{\bf 5.3.1}  The author has not yet coded 
{\tt AllConnectedSets} and {\tt AllMoConvexSets} which both feed on $CLP(G)$ (as opposed to $Geo(G)$).  It is to be expected that {\tt AllConnectedSets} in spe would stand up to {\tt BooleanConvert} as good as {\tt AllMetricSets} did. 

 I postponed {\tt AllConnectedSets} because from an application point of view (keyword: community detection) it is the most relevant of the four algorithms discussed in Section 5. It should therefore be pitted, on a common platform, against the
several algorithms for enumerating $Conn(G)$ that have recently been proposed in the literature (see [W2]).
 These "several algorithms" output the members of $Conn(G)$ one-by-one and predictably are inferior for real-life sized graphs $G$ with trillions\footnote{Even when compressed, dealing with trillions of subgraphs is usually not what is desired. See Subection 3.2  for ideas how to mitigate that.} of connected subsets. Help for setting up mentioned platform is welcome.

\section{Induced subgraphs that are cycle-free \\ or triangle-free or chordal or bipartite }

 By definition, a graph is {\it cycle-free} if it has no cycles at all. These graphs coincide with forests; in particular "cycle-free and connected"  means "tree". Yet we keep on using "cycle-free" to emphasize our agenda of investigating not only cycle-free subgraphs $G[X]$, but also three types of $G[X]$ which by definition are free of {\it  particular cycles}.
 
 Pleasantly in all of this it suffices to look at {\bf chordless} cycles. For starters, it is easily seen that every graph with a cycle also has a chordless cycle.
Given $G=(V,E)$ it therefore holds that $X \subseteq V$ induces a cycle-free subgraph $G[X]$ iff $C \not\subseteq X$ for all   chordless cycles  $C$.  Hence all cycle-free subgraphs $G[X]$ can be generated by applying the noncover $n$-algorithm to 
$CLC(G)$ (how to calculate $CLC(G)$ was discussed in Section 4). Upon invoking Theorem $1'$ we therefore derive

\vspace{5mm}
{\bf Theorem 3: }{\it Suppose the $h\ge 0$  chordless cycles of the graph $G=(V,E)$ are known. Then the $N$ sets $X\s V$ for which $G[X]$ is a forest, can be enumerated in time $O(Nh^2|V|^2)$.}

\vspace{2mm}
In the remainder we generalize  $G[X]$ from being a forest to being chordal (6.1), being bipartite (6.2), or being triangle-free (6.3). All that needs to be done is to feed particular subfamilies of $CLC(G)$ to the noncover $n$-algorithm.

\vspace{5mm}

{\bf 6.1} A graph is {\it chordal} [BM,p.235] iff each {\it long} cycle (i.e. of length $\ge 4$) has a chord. Hence all chordal subgraphs $G[X]$ can be generated by applying the noncover $n$-algorithm to the set $ LongCLC(G)$ of all long chordless cycles. In tandem with Theorem 3 we hence obtain

\vspace{3mm}

{\bf Theorem 4: }{\it Suppose the $h\ge 0$ {\bf long} chordless cycles of the graph $G=(V,E)$ are known. Then the $N$ sets $X\s V$ for which $G[X]$ is chordal, can be enumerated in time $O(Nh^2|V|^2)$.}

\vspace{3mm}

{\bf 6.1.1} Consider the toy graph $G_3=(V^3,E_3)$ in Figure 3A. Edge-wise it has the six cl cycles $\Gamma_1,...,\Gamma_6\s E_3$ listed in (9).

\begin{center}
\includegraphics[scale=0.5]{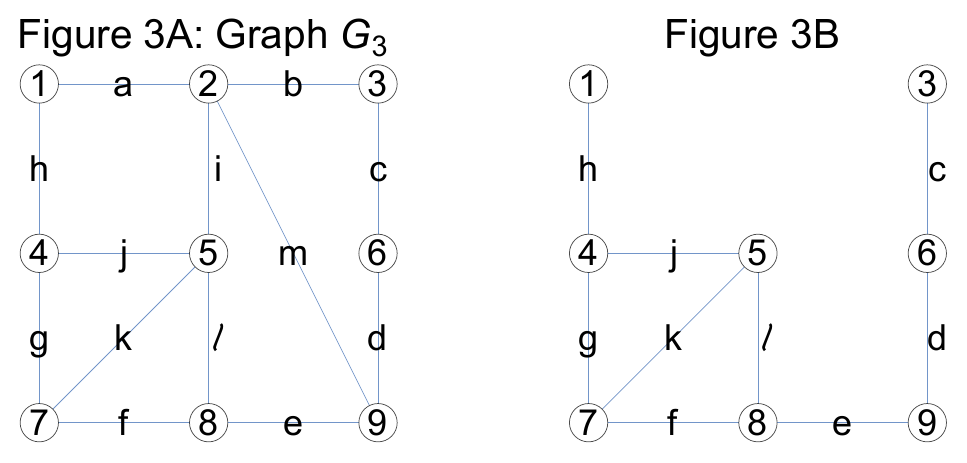}
\end{center}

\noindent
But of course we are concerned with the underlying {\it vertex} sets $C_1,...,C_6\s V$. Observe that $C_1,..,C_4$ are long chordless cycles, whereas $C_5,C_6$ are triangles. Moreover $C_1,...,C_4$ coincide with $A_1,...,A_4$ in (1).

\vspace{3mm}
(9) \quad $\begin{array}{llllll}
C_1 & = & \{1, 2, 4, 5\} & \quad \Gamma_1 & =& \{a, i, j, h\}\\
\\
C_2 & =& \{1, 2, 4, 7, 8, 9\} & \quad \Gamma_2 & =& \{a, m, e, f, g, h\}\\
\\
C_3 & =& \{2, 5, 8, 9\} & \quad \Gamma_3 & =& \{m, e, \ell, i\} \\
\\
C_4 & = & \{2, 3, 6, 9\} & \quad \Gamma_4 & = & \{b, c, d, m\}\\
\\
C_5 & =& \{4, 5, 7\} & \quad \Gamma_5 & =& \{j, k, g\}\\
\\
C_6 & = & \{5, 7, 8\} & \quad \Gamma_6 & =& \{\ell, f, k\} \end{array}$

\vspace{3mm}
\noindent

By Theorem 4, if $r_5,...,r_9$ are as in Table 1 then
the set $Chordal(G_3)$ of all chordal induced subgraphs of $G_3$ has the compressed representation $r_5\uplus r_6\uplus r_7\uplus r_8\uplus
 r_9$. In particular $|Chordal(G_3)|=431$ by (2). Observe that $r_5={\cal P}(Z)$ for $Z:=V\setminus\{2\}=\{1,3,4,5,6,7,8,9\}$. Hence $G_3[Z]$ (rendered in Fig. 3B) is the maximum size induced chordal subgraph.

 All  forests are chordal graphs. Thus  the induced subforests of $G_3$ are among the 431 induced chordal subgraphs. The final rows $r_8,r_9$  happen to
consist entirely of $|r_8|+|r_9|=25$ forests (see Table 1). How many of the remaining $431-25$ chordal subgraphs $G_3[X]$ are also forests? Rather than checking them one-by-one, one can
 impose (try) the two triangles $C_5=\{4,5,7\}$ and $C_6=\{5,7,8\}$ upon $r_5,r_6,r_7$.

\vspace{4mm}
{\bf 6.2} Recall that each graph with a cycle has a chordless cycle.  The following is slightly more subtle.

\begin{itemize}
    \item[(10)] {\it If a graph has no {\bf odd} chordless cycles, then it has no odd cycles at all.}
\end{itemize}

\noindent
To see this, let $C$ be any cycle with an odd number $n$ of edges. It suffices to prove the existence of an odd chordless cycle. There is nothing to show if $C$ is chordless itself. Otherwise pick any chord of $C$. It yields two obvious smaller cycles with $m+1$ and $(n-m)+1$ edges respectively. Because $(m+1)+(n-m+1)=n+2$ is odd, exactly one of the small cycles is odd. Iterating this we end up with a {\it chordless} odd cycle, possibly a triangle.

\vspace{3mm}
{\bf Theorem 5: }{\it Suppose the $h\ge 0$ {\bf odd} chordless cycles of the graph $G=(V,E)$ are known. Then the $N$ sets $X\s V$ for which $G[X]$ is bipartite, can be enumerated in time $O(Nh^2|V|^2)$.}

\vspace{3mm}
{\it Proof.} Recall that a graph is bipartite iff it has no odd cycles. According to (10) this is equivalent to having no {\it chordless} odd cycles.
Hence Theorem 5 once more follows from Theorem $1'$. $\square$

\vspace{4mm}
{\bf 6.3} A graph $G=(V,E)$ is {\it triangle-free} if it has no triangles. Loosely speaking triangle-free graphs constitute the antipode to chordal graphs because there {\it every} chordless cycle is a triangle. The triangle-free graphs not only comprise the bipartite graphs, but are generally well researched. To quote just one result, recall that each planar graph is 4-colorable; if additionally it is triangle-free then it is 3-colorable.
Recall from Section 4 that $Triangles(G)\s CLC(G)$ is the family of all triangles.
Feeding  $Triangles(G)$ to the noncover $n$-algorithm yields all triangle-free induced subgraphs. Accordingly it follows from Theorem $1'$:

\vspace{3mm}
{\bf Theorem 6: }{\it Suppose the $h\ge 0$ {\bf triangles} of the graph $G=(V,E)$ are known. Then the $N$ sets $X\s V$ for which $G[X]$ is triangle-free, can be enumerated in time $O(Nh^2|V|^2)$.}

\vspace{2mm}

As will be seen in Section 10.1, not only can $Triangles(G)$ be calculated fast, also the noncover $n$-algoritm can be trimmed considerably when its input is $Triangles(G)$.
Consequently, if e.g. $G$ in Theorem 6  is perfect (and thus has no long clcs, see [S,chapter 65]), then all  bipartite
($\stackrel{here}{=}$ triangle-free) subgraphs $G[X]$ can be calculated fast.

\section{Induced subgraphs that are chordless }

We call a graph {\it chordless} (not standard terminology) if all its connected components are either chordless paths or chordless cycles (in the standard sense). We stress that chordless paths of length 0 ( = isolated vertices) are also admitted as connected components.

On their own chordless graphs are rather boring (and they do not constitute a subclass of any of the four classes in Section 6). However, given $G=(V,E)$, it may be interesting to know all  chordless graphs of type $G[X]$.
For instance, if $G_1$ is from Figure 1 and $X_1,X_2,X_3\s V^1$ are defined by
$$(11)\quad X_1:=V^1\setminus\{2,7,9\},\ X_2:=V^1\setminus\{6,8\},\ X_3:=\{2,3,6,7,8,10\},$$
then $G_1[X_1],\ G_1[X_2],\ G_1[X_3]$ are the chordless subgraphs in Figures 1a, 1b, 1c whose edge-sets are rendered boldface.

 \vspace{3mm}
 {\bf 7.1}  For arbitrary $G=(V,E)$  and  $v\in V$ the {\it neighborhood} $NH(v)$ is the set of all vertices adjacent to $v$. For any fixed  $v\in V$ we say that $X\s V$ {\it satisfies  the $\{v\}$-constraint} if
$$(12)\quad v\in X\ \Ra\ |X\cap NH(v)|\le 2.$$
(In particular, the $\{v\}$-constraint is satisfied by $X$ if $v\not\in X$.) For all $X\s V$ we claim that: 
$$(13)\quad G[X]\ is\ chordless\ \LRa\  X\ satisfies\ all\ \{v\}\!-\!constraints\ (v\in V)$$
This follows from the  observation that each connected  graph with all vertices of degree $\le 2$ is either a chordless (!) path (possibly of length 0) or a chordless cycle.

\vspace{2mm}

\begin{tabular}{l|c|c|c|c|c|c|c|c|c|c|c|l} 
& 1 & 2 & 3 & 4 & 5 & 6 & 7 & 8 & 9 & 10 & 11& \\ \hline 
& & & & & & & & & & & &\\ \hline 
$r_1=$ & $a$ & $\epsilon(2)$ & $\epsilon(2)$ &$a'$  &$\epsilon'(2)$ & $\epsilon(2)$ &$\epsilon'(2)$  & $\epsilon'(2)$ &  &  & & p.7 \\ \hline 
& & & & & & & & & & & &\\ \hline 
$r_2=$ & $a$ & $\epsilon(2)$ & $\epsilon(2)$ &$2$  &$2$ & $\epsilon(2)$ &{\bf 0}  & $2$ &  &  & & p.2 \\ \hline 
$r_3=$ & $a$ & $\epsilon(2)$ & $\epsilon(2)$ &$a'$  &$\epsilon'$ & $\epsilon(2)$ &{\bf 1}  & $\epsilon'$ &  &  & & p.7 \\ \hline 
& & & & & & & & & & & &\\ \hline 
$r_4=$ & $2$ & {\bf 0} & $2$ &  & & $2$ &$0$  &  &  &  & & p.9 \\ \hline 
$r_5=$ & $a$ & {\bf 1} & $\epsilon$ &  & & $\epsilon$ &$0$  &  &  &  & & p.2 \\ \hline 
$r_3=$ & $a$ & $\epsilon(2)$ & $\epsilon(2)$ &$a'$  &$\epsilon'$ & $\epsilon(2)$ &$1$  & $\epsilon'$ &  &  & & p.7 \\ \hline 
& & & & & & & & & & & &\\ \hline 
$r_6=$ &  & $0$ &  &  & &  &$0$  &  & {\bf 0} &  & & final \\ \hline 
$(r_7'=$ &  & $0$ &  &  & &  &$0$  &  & {\bf 1} &  & & p.9) \\ \hline 
& & & & & & & & & & & &\\ \hline 
$r_7=$ &  & $0$ &  &  & &  &$0$  & $\epsilon(2)$ & $1$ & $\epsilon(2)$ & $\epsilon(2)$ & p.8 \\ \hline
$r_5=$ & $a$ & $1$ & $\epsilon$ &  & & $\epsilon$ &$0$  &  &  &  & & p.2 \\ \hline 
$r_3=$ & $a$ & $\epsilon(2)$ & $\epsilon(2)$ &$a'$  &$\epsilon'$ & $\epsilon(2)$ &$1$  & $\epsilon'$ &  &  & & p.7 \\ \hline 
\end{tabular}

\vspace{2mm}
{\sl Table 4: Computing $Chordless(G)$ with the $(a,\epsilon)$-algorithm}

\vspace{4mm}

Recalling the $(a,\epsilon(2))$-wildcard from Section 3  suppose that $G=(\{1,2,..\},E)$ is such that $NH(1)=\{2,3,4,5\}$. Then $(a,\epsilon(2),\epsilon(2),\epsilon(2),\epsilon(2),2,2,...,2)$ is the family of all $X\s V$ that satisfy the $\{1\}$-constraint. We therefore strive to impose these type of wildcards one-by-one akin to the $n$-wildcard in  Table 1.

Let us illustrate the details on $G_1$ when the order\footnote{Here "imposing a vertex $v$" is shorthand for imposing the $\{v\}$-constraint. The order of imposition is irrelevant (but influences the speed of the $(a,\epsilon)$-algorithm). The particular ordering is chosen to trigger certain illuminating effects.} of vertices to be imposed kicks off as $1,4,7,2,9,8,...$. For the sequel see Table 4, where for better visualization some of the don't-care 2's are replaced by blanks.

For starters, clearly each $X\in r_1$ satisfies both the $\{1\}$-constraint and the $\{4\}$-constraint. In order to impose $v=7$ on $r_1$, we split $r_1=r_2\uplus r_3$ as shown. Then vertex $2$ is pending in $r_2$ and $7$  remains pending in $r_3$. Upon splitting $r_2=r_4\uplus r_5$ as shown, 9 is pending in $r_4$ and $2$ remains pending in $r_5$. Upon splitting $r_4=r_6\uplus r_7'$  we find that $r_6$ happens to be final, i.e. all $X\in r_6$ yield chordless graphs $G_1[X]$. In fact $X_1$ in (11) is the largest member of $r_6$. What about $r_7'$? It turns out that we can impose its pending vertex 9 immediately (i.e. without splitting rows). The result is $r_7$, whose pending vertex is $8$.
And so it goes on. (The reader is invited to go just one step further: Impose vertex $8$ upon $r_7$.)

\vspace{3mm}

{\bf Theorem 7: }{\it Given $G=(V,E)$, the $N$ sets $X\s V$ for which $G[X]$ is chordless can be calculated in time $O(N|V|^4)$.}

\vspace{2mm}

{\it Proof.} By (13) the problem reduces to the imposition of all $\{v\}$-constraints $(v\in V)$. This amounts to imposing at\footnote{All vertices $v_0\in V$ with $deg(v_0)\le 2$ are harmless in the sense that {\it all} $X\s V$ trivially satisfy such $\{v_0\}$-constraints.} most $h:=|V|$ wildcards of type $(a,\epsilon(2),...,\epsilon(2))$. Since by Section 3 such wildcards can be reduced to $n$-wildcards, the $O(Nh^2|V|^2)$ of Theorem $1'$ applies. But here $O(Nh^2|V|^2)=O(N|V|^4)$. 
$\square$

\vspace{2mm}
It can happen  for wildcards based on $a,\epsilon,\epsilon(2)$ (and for other wildcards in Section 3) that many of the final rows end up being "shallow" 012-rows. But even so, these 012-rows may compress better than in a scenario where one opts for 012-rows from the start. As to "opt for", in the present scenario it is even {\it conceptually} simpler to embrace wildcards,
isn't it?

\vspace{2mm}
{\bf 7.2} Given $G=(V,E)$, let $Chordless(G)$ be the family of all sets $X\s V$ for which $G[X]$ is chordless. One verifies at once that $Chordless(G)$ is a set-ideal of ${\cal P}(V)$. Therefore, apart from representing $Chordless(G)$ with wildcards one could capture it by determining the maximal members of this set-ideal.
While for many set-ideals in combinatorics calculating their maximal members (the so called {\it facets}) is well researched, this remains an open question for $Chordless(G)$.

\vspace{2mm}
{\bf 7.3}  What about switching $\le 2$ in (12) for $\le 1$? Then one gets (still in polynomial total time) all matchings of $G=(V,E)$ which are "isolating" in the sense that  distinct edges in the matching are never joined by an edge of $E$.
For instance, the matching $\{\{1,2\},\{3,6\},\{4,5\},\{8,9\}, \{10,11\}\}$ of $G_1$ is not isolating, but $\{\{3,6\},\{4,5\},\{10,11\}\}$ is.

What about switching $\le 2$ in (12) for $=2$? Then one only gets chordless graphs $G[X]$ all of whose connected components are (isolating) chordless cycles.  While this may be desirable for some applications, the prize is that $a,\epsilon,\epsilon(2)$ give way to $g,g(2)$ (see Section 3), and the polynomial total time in Theorem 7 is lost.

\vspace{2mm}
{\bf 7.4} A subset $X\s V$ is an {\it anticlique} (or independent set) of a graph $G=(V,E)$ if no two vertices in $X$ are adjacent to each other. Article [W4] is dedicated to calculating (in compressed fashion) the family $Acl(G)$ of all anticliques of $G$. The basic idea is to use the $(a,c)$-wildcard mentioned in (4) of Section 3. Namely, each $v\in V$ gives rise to a wildcard $(a,c,c,\ldots,c)$ where $a$ matches $v$ and $(c,c,...,c)$ matches $NH(v)$. Imposing these wildcards one after the other, akin to Table 4, one obtains $Acl(G)$. In view of Section 3 one can view $Acl(G)$ as the modelset of some Horn CNF.

Anticliques will  come up a few times in the present article. To begin with, for $X\in Chordless(G)$ let $\{X_1,...,X_s\}$ be the connected components of $G[X]$. It then follows that
$$(14)\quad Acl(G[X])=\big\{Y_1\uplus\cdots\uplus Y_s:\ (\forall 1\le i\le s)\ Y_i\in Acl(G[X_i])\big\},$$
i.e. $Acl(G[X])$ reduces nicely to the smaller parts  $Acl(G[X_i])$. What is more,
these parts have a succinct structure. For instance, for a chordless paths $X_i$ of length $m$ one has $|Acl([G[X_i])|=Fib(m+3)$, where $Fib(1):=0,Fib(2)=1,Fib(n):=Fib(n-1)+Fib(n-2)\ (n\ge 3)$ are the classic Fibonacci numbers.

\section{Packing cliques}

  Let $G=(V,E)$ be a graph. A partition $\Pi=\{V_1,...,V_k\}$ of $V$ is a {\it packing of type T subgraphs} if for each {\it proper}  ( := non-singleton) part $V_i\in P$ it holds that the induced subgraph $G[V_i]$ is of type $T$.
In this Section the packings $\Pi=\{V_1,...,V_k\}$ of type $T$ graphs are {\it clique-packings (Cli-Pac)}. Thus $\Pi$ is a partition of the vertex-set of $G=(V,E)$ and each $V_i$ is a clique of $G$. By definition $E(\Pi)$ is the set of all edges that belong to one of the cliques $V_i\in \Pi$.
Consequently the connected components of the graph $G^\Pi:=(V,E(\Pi))$ are the sets $V_i\in\Pi$.
Each  $V_i$ is a maximal clique wrt $G^\Pi$ but not necessarily maximal wrt $G$.

\begin{center}
\includegraphics[scale=0.9]{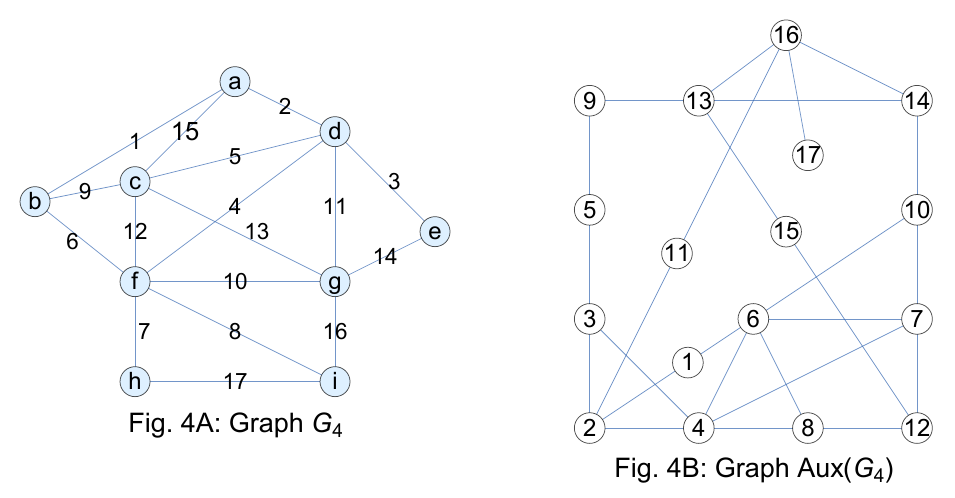}
\end{center}

For instance $G_4$ in Fig. 4A  has the   Cli-Pacs $\Pi_1=\{\{c,d,f,g\},\{a,b\},\{h,i\},\{e\}\}$ 
(with $E(\Pi_1)=\{4,5,10,11,12,13,1,17 \}$) and $\Pi_2=\{\{a,b,c\},\{d,e,g\},\{f,h,i\}\}$ \\ 
(with $E(\Pi_2)=\{1,9,15,3,11,14,7,8,17 \}$).

\vspace{2mm}
{\bf 8.1} At this point the reader may wonder: Can the chordless graphs of Section 7 also be considered as a kind of packing?
Before we continue with Cli-Pacs, let us clarify this.
Define the property $T$ as "being a chordless path or chordless cycle". Take any $X\in Chordless(G)$ and let $\Pi_X$ be the partition whose parts are the connected components of $G[X]$. Then $\Pi_X$ is indeed a packing of type $T$ graphs. This packing has the additional feature of being {\it isolating} in the sense that there are  no edges connecting different parts of $\Pi_X$; see 7.3 for a concrete example. This is in stark contrast to e.g. $\Pi_1$ where there are {\it several} edges between the parts $\{c,d,f,g\}$ and $\{a,b\}$. 

Whatever the type $T$, each isolating packing of type $T$ subgraphs of $G$ can always be reduced to the consideration  of induced subgraphs $G[X]$. The packings dealt with in Sections 8 and 9 are not isolating.

\vspace{2mm}
{\bf 8.2} Back to Cli-Pacs.  Two types of constraints will be crucial. First,
coupled to each {\it edge-triangle} $\{\alpha,\beta,\gamma\}\s E$ (definition clear) we introduce a certain "constraint" $[\{\alpha,\beta,\gamma\}]$. By definition $Y\s E$ {\it satisfies}
$[\{\alpha,\beta,\gamma\}]$ if
$$(15)\quad |Y\cap \{\alpha,\beta,\gamma\}|\ge 2\ \Ra\ \{\alpha,\beta,\gamma\}\s Y$$
In plain language, either $Y$ contains at most one of $\alpha,\beta,\gamma$, or all of them. Second, whenever $\{\alpha,\beta\}\s E$ is the edge set of a {\it chordless} path, it triggers another type of constraint $[\{\alpha,\beta\}]$. By definition $Y\s E$ {\it satisfies} $[\{\alpha,\beta\}]$ if
$$(16)\quad |Y\cap \{\alpha,\beta\}|\le 1.$$

\noindent
Let $\Sigma_G$ be the family of all these constraints (both types). 
A subset $Y\s E$ is $\Sigma_G$-{\it closed} if $Y$ satisfies all constraints in $\Sigma_G$.

\vspace{5mm}
{\bf Lemma 8: }
\begin{itemize}
\item[(a)] {\it If $\Pi$ is a Cli-Pac of $G=(V,E)$ then $E(\Pi)$ is $\Sigma_G$-closed. }
\item[(b)] {\it If $Y\s E$ is
any $\Sigma_G$-closed set then the connected components of $G':=(V,Y)$ yield a Cli-Pac $\Pi$ of $G$ with $E(\Pi)=Y$.}
\end{itemize}
\vspace{3mm}

{\it Proof.} (a) Fix any type 1 constraint $[T]:=[\{\alpha,\beta,\gamma\}]$ in $\Sigma_G$. Assuming that (say) $\{\alpha,\gamma\}\s E(\Pi)$ we must show that $T\s E(\Pi)$. Let $a,b,c$ be the vertices of the edge triangle $T$. From $\{\alpha,\gamma\}\s E(\Pi)$ follows that
$a,b,c$ are in the same connected component of $G^\Pi$, i.e. the same clique $V_i$ of $G$. Hence $T\s E(\Pi)$. Likewise fix any type 2 constraint $[\{\alpha,\beta\}]$ in $\Sigma_G$.
Let $(a,b,c)$ be the coupled chordless path. Hence  $\alpha=\{a,b\}$ and $\beta=\{b,c\}$ and $\{a,c\}\not\in E$. We need to show that
$\{\alpha,\beta\}\not\s E(\Pi)$. If not, then $a,b,c$ were again in the same clique $V_j$ of $G$, and in particular $\{a,c\}\in E(\Pi)$. This contradicts $\{a,c\}\not\in E$.

\vspace{2mm}
(b) Let $V_0\s V$ be any fixed connected component of $G'=(V,Y)$. There is nothing to show if $|V_0|=1$, and so let $V_0$ be proper. We need to show that $V_0$ is a clique of $G$. Let $\{v_1,...,v_k\}$ be a maximal clique of $G'$ contained in $V_0$; here  $k>1$ since $V_0$ is proper. It suffices to show $\{v_1,...,v_k\}=V_0$. Thus let us derive a contradiction from assuming the existence of $v\in V_0\setminus\{v_1...,v_k\}$. Since $V_0$ is connected there is a path in $G'$ from $v$ to $v_k$. Upon relabeling we can assume that $v$ is adjacent to $v_k$.
 Because  $\alpha:=\{v_{k-1},v_k\}$ {\it and} $\beta:=\{v_k,v\}$ belong to $Y$, and $Y$ is $\Sigma_G$-closed by assumption, the path $(v_{k-1},v_k,v)$ must have\footnote{Put another way, if there was no chord then $[\{\alpha,\beta\}]\in\Sigma_G$. The $\Sigma_G$-closedness of $Y$ then yields the contradiction ($\alpha\not\in Y$ or $\beta\not\in Y$).} the chord $\{v_{k-1},v\}\in Y$. Likewise it follows from $\{v_{k-2},v_{k-1}\},\{v_{k-1},v\}\in Y$ that $\{v_{k-2},v\}\in Y$, and so forth.
 Therefore $\{v_1,...,v_k,v\}$ is a clique, which contradicts the maximality of $\{v_1,...,v_k\}$. $\square$

\vspace{2mm}
Observe that both constraint types (15) and (16)  are about triangles. One involves the vertices, the other the edges of the triangle.

\vspace{2mm}
{\bf Theorem 9: }{\it  The $N$ many Cli-Pacs of a graph $G=(V,E)$ can be calculated (in compressed format) in time $O(N|V|^{10})$.}

\vspace{2mm}
{\it Proof.} According to Lemma 8 the $N$ Cli-Pacs bijectively match the $\Sigma_G$-closed sets $Y\s E$. Clearly $Y\s E$ satisfies the type 2 constraint $[\{\alpha,\beta\}]$ iff $Y$ is an $\{\alpha,\beta\}$-noncover. As to type 1 constraints, $Y\s E$ satisfies
$[\{\alpha,\beta,\gamma\}]$ iff $Y$ is closed wrt the three implications $\{\alpha,\beta\}\ra\{\gamma\},\ \{\alpha,\gamma\}\ra\{\beta\},\ \{\beta,\gamma\}\ra\{\alpha\}$. We see that being $\Sigma_G$-closed amounts to be a model of some Horn formula $f$ with a variable set $W$ equicardinal to $E$, and which has at most $h:=\binom{|V|}{2} + 3\binom{|V|}{3}$ clauses. According to Theorem 1 one can enumerate the $N$ models of $f$ in time $O(Nh^2|W|^2)=O\big(N\cdot (|V|^3)^2\cdot(|V|^2)^2\big)=O(N|V|^{10})$. $\square$

\vspace{2mm}
{\bf 8.2.1.} We mention that $|E(\Pi_2)|=9$ maximizes $|E(\Pi)|$ when $\Pi$ ranges over all Cli-Pacs of $G_4$. 
As to generally producing only those Cli-Pacs $\Pi$ with $E(\Pi)$ sufficiently large, say $|E(\Pi)|\ge k$,  since
Cli-Pacs are based on Horn-formulas, the technique of 3.2 allows to generate all such Cli-Pacs in polynomial total time.

\vspace{2mm}
{\bf 8.3} By Theorem 9 the Cli-Pacs of a graph $G$ can be viewed as the $\Sigma_G$-closed sets of some implication-family $\Sigma_G$, and thus as the models of some Horn formula $f_G$. What about a Boolean formula $\ol{f}_G$ whose models match the {\it individual} cliques of $G$? In order to succeed consider the {\it complementary} graph $\ol{G}:=(V,\ol{E})$ which is defined by $\ol{E}:=\binom{V}{2}\setminus E$ where $\binom{V}{2}$ is the set of all 2-element subsets of $V$. Clearly it holds for each $X\s V$ that:
$$(17)\quad X\ is\ a\ clique\ of\ G\ \LRa\ X\ is\ an\ anticlique\ of\ \ol{G}. $$
\noindent
As seen in 7.4, there is some Horn CNF $F_{\ol{G}}$ whose models match the anticliques of $\ol{G}$. Hence $\ol{f}_G:=F_{\ol{G}}$ does the job, in view of (17).

It also follows from (17) that all anticlique-packings of $G$ coincide with the Cli-Pacs of $\ol{G}$, which can be calculated as described above. As is well known, anticlique-packings of $G=(V,E)$ and (proper) colorings $c:V\to\{1,2,..,k\}$ are the same thing. Namely, each $k$-coloring $c$ yields the anticlique-packing
$\Pi_c:=\{c^{-1}(1), c^{-1}(2),...,c^{-1}(k)\}$; and conversely each anticlique-packing $\Pi$ yields and obvious $k$-coloring $c_{\Pi}:V\to\{1,2,...,k\}$. The author doubts that this observation is relevant for calculating the color number $\chi(G)$ of $G$ (but it might be if by whatever reason {\it all} $\chi(G)$-colorings of $G$ are required).

\section{Packing connected sets}

Given $G=(V,E)$, recall that $X\s V$ being connected means that $G[X]$ is connected. In Section 5 we enumerated all connected sets $X$ individually. In Section 9 we {\it pack} them and,
surprisingly, will fare better what concerns polynomial total time!
Having glimpsed  at naive attempts in 9.1, Subsection 9.2 embraces more efficient techniques to generate all packings of connected sets (=: {\it Conn-Pacs}).

\vspace{4mm}
{\bf 9.1} Consider $G_5=(V^5,E_5)$ in Fig. 5.1. Let us pick at random any partition of $V:=V^5$, such as
$\{\{1,2,5,6,11,12,15,16\},\{3,4,7,8,9,10,13,14\}\}$. Thus $V$ splits  into bright and dark vertices. Both the set of bright, and the set of dark vertices, evidently induce {\it disconnected} subgraphs. This is the expected outcome when trying to find Conn-Pacs by trial and error.

 One may hence be led to proceed differently: take all individual connected sets $X$ (obtained in whatever way) and try to combine them to Conn-Pacs. However, this is hopeless as well\footnote{It is like taking the billion  pieces of a million different puzzle games, mixing them thoroughly, and then attempting to solve {\it all} puzzle games.}. 

\begin{center}
\includegraphics[scale=0.87]{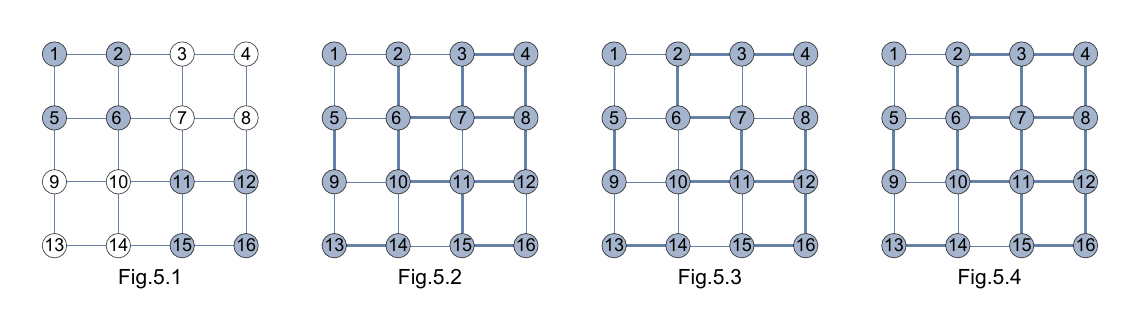}
\end{center}

{\bf 9.1.1} The following somewhat improves upon these simple-minded attempts. Take any set of edges $K\s E_5$ and let $V_1,...,V_t$ be the connected components of the graph $(V,K)$. If $K$ is the set of bold edges in Fig.5.2, then the induced Conn-Pac is $\Pi_0=\{V_1,...,V_t\}:=(\{1\},\{5,9\},\{13,14\},\{2,3,4,6,7,8,10,11,12,15,16\})$. It is clear that {\it every} Conn-Pac  can be obtained this way. 

Trouble is, it may be obtainable in multiple ways. For instance, the set $K'\s E_5$ of bold edges in Fig.5.3 also yields $\Pi_0$.
The good news is, among all edge-sets inducing a fixed Conn-Pac $\Pi$ of a graph $G=(V,E)$, there is a largest edge-set, call it $E(\Pi)$. For instance $E(\Pi_0)$ is the set of bold edges in Fig.5.4. Generally

$$(18)\quad E(\Pi):=\Bigl\{\{u,v\}\in E:\ (\exists i\le t)\ u,v\in V_i \Bigr\} $$

\vspace{3mm}
{\bf 9.2} Recall that $CLC(G)$ is the set of all chordless cycles of $G=(V,E)$. Here comes a seemingly clumsy way to characterize edge-sets of type (18), yet it will link the whole matter to implications $A\ra B$ (as defined in Section 2).

\vspace{3mm}
{\bf Lemma 10: }
{\it Take any $K\s E$. It is of type $K=E(\Pi)$ for some Conn-Pac $\Pi$ iff for each $C\in CLC(G)$  it holds that $|C\setminus K|\neq 1$.}

\vspace{2mm}
\noindent Before embarking on the proof, fix any $K\s E$ and any  $C\in CLC(G)$. One readily checks that these two conditions are equivalent:

\begin{itemize}
\item[(19)] $|C\setminus K|\neq 1$
    \item[(20)]{\it If $C$  gets "nearly" swallowed  by $K$ (meaning that $(C\setminus K)\s\{\alpha\}$), then 
    $C$ gets "wholly" swallowed by $K$ (so $\alpha\in C$, i.e. $C\s K$).}
\end{itemize}

{\it Proof of Lemma 10.} As to $\Ra$, consider $E(\Pi)$ for some Conn-Pac $\Pi=\{V_1,...,V_t\}$. Fix any $C\in CLC(G)$. In order to verify (20) for $C$ and $K:=E(\Pi)$, assume that \\
 $C=\{\{u_1,u_2\},..,\{u_{s-1},u_s\},\{u_s,u_1\}\}$ and that 
$\{\{u_1,u_2\},...,\{u_{s-1},u_s\}\}\s E(\Pi)$.  We must show that $\{u_{s},u_1\}\in E(\Pi)$.
Since there is a path from $u_1$ to $u_s$ whose edges lie in $E(\Pi)$, the vertices $u_1,u_s$ lie in the same connected component of $(V,E(\Pi))$, i.e. $u_1,u_s\in V_i$ for some $V_i\in \Pi$. From $\{u_1,u_s\}\in E$ and definition (18) follows  $\{u_s,u_1\}\in E(\Pi)$. 

\vspace{3mm}
As to $\Leftarrow$, suppose that $K\s E$ is such that (20) holds for each $C\in CLC(G)$. Let $\Pi_0:=\{V_1,...,V_t\}$ be the Conn-Pac whose parts $V_i$ are the connected components of the graph $G'=(V,K)$. We claim that $K=E(\Pi_0)$. 
Fix any edge $\{u_1,v\}\in E$ such that $u_1,v$ are in the same component of $\Pi_0$, say $u_1,v\in V_3$. If we can show that $\{u_1,v\}\in K$, then $K=E(\Pi_0)$ as claimed. To begin with, by definition of $\Pi_0$ there is a path $P$ from $u_1$ to $v$  whose edge-set lies in $K$, say 
$E(P)=\{\{u_1,u_2\},\{u_2,u_3\},...,\{u_s,v\}\}\s K$. From $\{v,u_1\}\in E$ follows that $E(P):=P\cup\{\{v,u_1\}\}$ is a cycle. If $C'$ happens to be chordless  then it follows from (20)  that $\{v,u_1\}\in K$. If $C'$ has chords then, according to (21) below, it still holds that $\{v,u_1\}\in K$. $\square$

\begin{itemize}
    \item[(21)]{\it Let $G=(V,E)$ a graph and $K\s E$ arbitrary. If (20) holds for all $C\in CLC(G)$, then
    (20) holds for {\bf all} cycles $C'$ of $G$.}
\end{itemize}

In order to prove (21) by way of contradiction, let $C'$ be a cycle of minimum cardinality that violates (20), i.e. it holds that $C\setminus K=\{\alpha\}$ for some $\alpha\in C$. Since $C'$ violates (20), it has some chord $\beta$. This chord splits $C'$ into two parts. Each part together with $\beta$ yields a cycle; call them $C_1'$ and $C_2'$. Evidently $C_1'\setminus\{\beta\}$ and  $C_2'\setminus\{\beta\}$ are contained in $C'$. Say $\alpha\in C_1'$, and so $\alpha\not\in C_2'$. From $C_2'\setminus\{\beta\}\s C'\setminus\{\alpha\}\s K$ and $|C_2'|<|C'|$ follows $\beta\in K$. This together with $(C_1'\setminus\{\beta\})\setminus\{\alpha\}\s C'\setminus\{\alpha\}\s K$ implies that $C_1'\setminus\{\alpha\}\s  K$. But $|C_1'|<|C'|$ implies $\alpha\in K$. This contradicts $C\setminus K=\{\alpha\}$ and thus proves (21).

\vspace{3mm}
 Defining the  implication-family

$$(22)\quad \Sigma^{G}:=\Bigl\{(C\setminus \{e\})\to\{e\}:\ C\in 
CLC(G),\ e\in C\Bigr\}$$

\noindent
it follows from Lemma 10 that the $\Sigma^G$-closed edge-sets $K\s E$ are exactly the edge-sets of type $E(\Pi)$.

\vspace{3mm}
{\bf Theorem 11: }{\it Suppose the $h\ge 0$ chordless  cycles of the graph $G=(V,E)$ are known. Then the $N$ Conn-Pacs of $G$ can be enumerated in time $O(Nh^2|V|^2|E|^2)$.}

\vspace{3mm}
{\it Proof.} By assumption $|CLC(G)|=h$, and so $h_1:=|\Sigma_G|\le h|V|$.
We saw that  the $N$  Conn-Pacs of $G$ bijectively match the $\Sigma^G$-closed subsets of $E$. According to Theorem $1''$ the latter can be enumerated in time $O(Nh_1^2w^2)=O(N(h|V|)^2|E|^2)=O(Nh^2|V|^2|E|^2)$. $\square$

\vspace{2mm}
Let $k>0$ be fixed. By the same reason as in 8.2.1 one can, without sacrificing polynomial total time, restrict the output in Theorem 11 to the Conn-Pacs $\Pi$ with $|E(\Pi)|\ge k$.

\vspace{3mm}
{\bf 9.3}  For readers familiar with matroids [S,chapter 39] we mention that the members of
$${\cal F}_E(G):=Cl(\Sigma^G)$$
are exactly the so-called {\it flats} of a graphic matroid\footnote{Akin to Thm.11 the flat lattice ${\cal F}(M)$ of {\it any} matroid $M$ can be enumerated [M], but by wholly different means (e.g. relying on a independence oracle). Assuming the independence oracle is polynomial time, also ${\cal F}(M)$  gets enumerated in polynomial time, albeit  one-by-one.} $M_G$ coupled to the graph $G=(V,E)$. Like every closure system, ${\cal F}_E(G)$  is a lattice [S,p.233] where the meet $F_1\wedge F_2$ of $F_1,F_2\in{\cal F}_E(G)$ is simply $F_1\cap F_2$. 
It turns out that
$${\cal F}_V(G):=\Big\{\Pi:\ \  \Pi\ is\ a\ Conn\!-\!Pac\ of\ G\Big\}$$
has a lattice structure as well.
In fact $F\mapsto \Pi(F)$ defines a lattice isomorphism from ${\cal F}_E(G)$ to ${\cal F}_V(G)$. Here by definition the parts of the partition $\Pi(F)$ of $V$ are the connected components of the graph $(V,F)$. The inverse lattice isomorphism is given by
$\Pi\mapsto E(\Pi)$.
One can view $Conn(G)$ as a subset of ${\cal F}_V(G)$. More precisely, each $X\in Conn(G)$ yields the member $\{\{X\}\}\cup\{\{y\}:\ y\in V\setminus X\}$ of ${\cal F}_V(G)$.

 Although the lattices ${\cal F}_V(G)$ and ${\cal F}_E(G)$ are isomorphic, meets in the former are more subtle than in the latter. To witness, $\Pi_1\wedge\Pi_2$ is the coarsest\footnote{For instance, the intersection of the two connected sets $X_1=\{1,2,3,4,5\}$ and $X_2=\{2,3,4,5,6\}$ of $G_2$ is disconnected (Fig. 2B). In contrast, if $\Pi_1:=\{X_1,\{6\}\}$  and $\Pi_2:=\{X_2,\{1\}\}$ are the coupled Conn-Pacs in (i.e. members of ${\cal F}_V(G_2)$), then $\Pi_1\wedge\Pi_2=\{\{1\},\{2,3\},\{4,5\},\{6\}\}$ is still a Conn-Pac! }  Conn-Pac which is  finer than both $\Pi_1$ and $\Pi_2$.

\vspace{3mm}
{\bf 9.4} Similar to 7.2, ${\cal F}_E(G)$ is not only  captured  by $\Sigma_G$), but also by the maximal members $H$ of ${\cal F}_E(G)\setminus\{E\}$. We like to call them the {\it edge-hyperplanes} $H$. Their images under the lattice isomorphism 
${\cal F}_E(G)\to {\cal F}_V(G)$ are the 
 {\it vertex-hyperplanes}, defined as those $\Pi\in {\cal F}_V(G)$ with\footnote{In particular, if $n:=|V|$, then there can be at most $\frac{2^{n}-2}{2}$ vertex-hyperplanes.} $|\Pi|=2$.
In brief, each vertex-hyperplane $\Pi$ yields the edge-hyperplane $E(\Pi)$, and each edge-hyperplane $H$ yields the vertex-hyperplane $\Pi(H)$.

One can obtain ${\cal F}_E(G)$  from the edge-hyperplanes by intersecting them in all possible ways; doing it more carefully, this can be done in linear total time\footnote{This holds for any closure system if its meet-irreducible elements are known. In our scenario the meet-irreducibles are the edge-hyperplanes.}.

The bad news is, this method outputs the flats one-by-one, and their may be zillions of them. Furthermore, how to obtain the edge-hyperplanes in the first place? Fortunately, this bit has a more uplifting answer. Namely, a {\it cutset} of a connected graph $G=(V,E)$ is any set $X\s E$ such that the graph $(V,E\setminus X)$ is disconnected. Furthermore, $X$ is {\it minimal} if each proper subset of $X$ fails to be a cutset. 
The  minimal cutsets can be generated in linear total time  according to [SA].

As is well known, the edge-hyperplanes $H\in {\cal F}_E(G)$ are exactly the sets $E\setminus X$  when $X$ ranges over all minimal cutsets of $G$.

 \vspace{4mm}
 {\bf 9.5} To hint at just one application, let $G=(V,E)$ be given. Here $G$ must be connected in order to have $\{V\}\in {\cal F}_V(G)$. Suppose a biologist comes to us ( = the mathematicians) with an {\it ordinary} partition $\pi_0=\{Y_1,...,Y_k\}$ of $V$. Most parts $Y_i$ may {\it not} be connected subsets of $G$. By whatever motivation the biologist seeks a Conn-Pac which is as "similar" to $\pi_0$ as possible. Pleasantly, among all Conn-Pacs in ${\cal F}_V(G)$ which are coarser than $\pi_0$ (such as $\{V\}$), there is indeed a finest (=most similar) one, call it $\Pi_0$.  Assuming that we {\it know} all vertex-hyperplanes, $\Pi_0$ can be calculated as follows. Say $\Pi_1,...,\Pi_t$  are those vertex-hyperplanes that happen to be coarser than $\pi_0$. It then holds that $\Pi_0=\Pi_1\wedge\cdots\wedge\Pi_t$.

 But how to obtain all vertex-hyperplanes in the first place? As described in 9.4, the {\it edge}-hyperplanes $H_1,H_2,...$ are obtained from the minimal cutsets of $G$, and so the vertex-hyperplanes originate by determining the connected components of  the graphs $(V,H_1),(V,H_2),...$.

\section{The noncover-algorithm and its variants: Part 2}

In Subsection 10.1 we show how the particular application of the  noncover $n$-algorithm that occurs in Theorem 6 can be trimmed; however this has not yet be programmed and does not improve the $O(Nh^2|V|^2)$ bound. In 10.2 we show how the particular application of the  Horn $n$-algorithm that occurs in Theorem 9 can be trimmed. In both 10.1. and 10.2 we exploit that applying the noncover $n$-algorithm to a family $\cal S$ of 2-element sets amounts to finding all anticliques of a graph. Finally  10.3 is dedicated to accelerate the calculation of all Conn-Pacs of a graph.

\begin{center}
\includegraphics[scale=0.9]{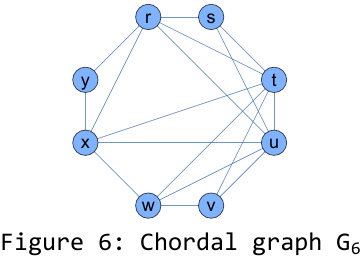}
\end{center}

\vspace{3mm}
{\bf 10.1} We start by generating all triangles of the  graph $G_6=(V^6,E_6)$ in Figure 6 (which is copied from [BM,p.235]).  Since $r$ is the lexicographic first letter in $V^6=\{r,s,..,x,y\}$ we look at $NH(r)=\{s,t,u,x,y\}$ and readily\footnote{Any two members of $NH(r)$ that happen to be adjacent yield a triangle.} find all six triangles that contain $r$. Next we evaluate $NH(s)$ wrt $G_6[V^6\setminus\{r\}]$, then $NH(t)$ wrt $G_6[V^6\setminus\{r,s\}]$, and so on. This yields the following groups of triangles:

\begin{itemize}
    \item[(23)] $\{r,s,t\},\ \{r,s,u\},\ \{r,t,u\},\ \{r,t,x\},\ \{r,u,x\},\ \{r,x,y\},$
    \item[] $\{s,t,u\},$
    \item[] $\{t,u,v\},\ \{t,u,w\},\ \{t,u,x\},\ \{t,v,w\},\ \{t,w,x\},$
    \item[] $\{u,v,w\},\ \{u,w,x\}$
\end{itemize}

{\bf 10.1.1} Consider row $\rho_1$ in Table 5. It is clear that each $Y\in \rho_1$ {\it avoids}\footnote{Here "avoids" is shorthand for "is a noncover of"} each triangle in (23) of type $\{r,\ast,\ast\}$. Because of $nnn$ in $\rho_1$, set $Y$ avoids $\{s,t,u\}$ as well. Hence the triangles pending to be imposed are the ones of type $\{t,\ast,\ast\}$. 

As to $\rho_2$, this is more complicated. To begin with, put $\rho_2':=(1,2,...,2)$. Evidently $Z\in \rho_2'$ avoids all triangles $\{r,\ast,\ast\}$ iff $Z\cap\{s,t,u,x,y\}$ is an anticlique  of the graph in Figure 7.1.  
The vertices of this graph  match the blanks in row $\rho_2$. Eventually these blanks will get filled by exactly those bitstrings that encode\footnote{Here we encode the anticliques ad hoc, see 10.2 on how to do this systematically.} the anticliques of this graph.
This is how the comment "pending Fig.7.1" on the right of $\rho_2$ is to be understood. The further comment "then $\{s,t,u\}$" conveys that afterwards, according to the  ordering in (23), the triangle $\{s,t,u\}$ is to be processed. 

As in Table 1 and 3 we keep on processing the top rows of the  stack. As to processing $\rho_1$, one has $\rho_1=\rho_{1,1}'\cup\rho_{1,2}':=(0,{\bf 2,0,2},2,2,2,2)\uplus (0,{\bf n,1,n},2,2,2,2)$. Because of $3\in zeros(\rho_{1,1}')$ all $y\in\rho_{1,1}'$ avoid the  type $\{t,*,*\}$ triangles. Shrinking $\rho_{1,1}'$ to $\rho_{1,1}$
(see Table 5) one achieves that additionally all $Y\in\rho_{1,1}$ avoid $\{u,v,w\}$, and so the last pending triangle for $\rho_{1,1}$ is $\{u,w,x\}$. Akin to above $Z\in \rho_{1,2}'$ avoids all triangles $\{t,*,*\}$ iff\footnote{The $\underline{n}$ in $\rho_{1,2}$ is underlined to indicate its obligations both  wrt $(n,\underline{n})$ and wrt Fig. 7.2.}
$Z\cap\{u,v,w,x\}$ is an anticlique of the graph in Fig. 7.2.

\vspace{4mm}
\begin{tabular}{l|c|c|c|c|c|c|c|c|c}
& r & s & t & u & v & w & x & y &\\ \hline
&  &  &  &  &  & & & & \\ \hline
${\cal P}(V^6)=$ & 2 & 2 & 2 & 2 & 2 & 2  & 2 & 2 & \\ \hline
 &  &  &  &  &  & & & & \\ \hline
$\rho_1=$ & {\bf 0} & $n$ & $n$ & $n$ & 2 & 2 & 2& 2 &
pending $\{t,\ast,\ast\}$ \\ \hline

$\rho_2=$ & {\bf 1} &  & & & 2 & 2&  &  & pending Fig.7.1, then $\{s,t,u\}$\\ \hline

&  &  &  &  &  & & & & \\ \hline
$\rho_{11}=$ & 0 & 2 & {\bf 0} & $n'$ & $n'$ & $n'$ & 2  & 2 & 
pending $\{u,w,x\}$ \\ \hline

$\rho_{12}=$ & 0 & $n$ & {\bf 1} & $\underline{n}$ &  & &  & 2 & 
pending Fig 7.2, then $\{u,\ast,\ast\}$ \\ \hline
$\rho_2=$ & {\bf 1} &  &  & & 2 & 2&  &  & pending Fig.7.1, then $\{s,t,u\}$\\ \hline

&  &  &  &  &  & & &  &\\ \hline
$\rho_{111}=$ & 0 & 2 & 0 & $n$ & {\bf 0} & $n$ & $n$  & 2 & final, card=28 \\ \hline
$\rho_{112}=$ & 0 & 2 & 0 & $n$ & {\bf 1} & $n$ & 2  & 2  & final, card=24 \\ \hline
$\rho_{12}=$ & 0 & $n$ & {\bf 1} & $\underline{n}$ &  & &  & 2 & pending Fig.7.2, then $\{u,\ast,\ast\}$ \\ \hline
$\rho_2=$ & {\bf 1} &  &  & & 2 & 2&  &  & pending Fig.7.1, then $\{s,t,u\}$\\ \hline

&  &  &  &  &  & & &  &\\ \hline
$\rho_{121}=$ & 0 & 2 & 1 & {\bf 0} &  & &  & 2 & pending Fig.7.3\\ \hline
$\rho_{122}=$ & 0 & 0 & 1 & {\bf 1} & 0 &0 &0  &2  & final, card=2 \\ \hline
$\rho_2=$ & {\bf 1} &  &  & & 2 & 2&  &  & pending Fig.7.1, then $\{s,t,u\}$\\ \hline

&  &  &  &  &  & & &  &\\ \hline
$\rho_{1211}=$ & 0 & 2 & 1 & 0 & 2 & {\bf 0}& 2 & 2 & final, card=16\\ \hline
$\rho_{1212}=$ & 0 & 2 & 1 & 0 & 0 &{\bf 1} & 0 & 2 & final, card=4\\ \hline
$\rho_2=$ & {\bf 1} &  &  & & 2 & 2&  &  & pending Fig.7.1, 
then $\{s,t,u\}$\\ \hline

&  &  &  &  &  & & &  &\\ \hline
$\rho_{21}=$ & 1 &  &  & & 2 & 2& {\bf 0} & 2 & pending Fig.7.4, 
then $\{t,\ast,\ast\}$ \\ \hline
$\rho_{22}=$ & 1 & 2 & 0 &0 & 2 & 2& {\bf 1} & 0 & final, card=8\\ \hline

&  &  &  &  &  & & &  &\\ \hline
$\rho_{211}=$ & 1 &{\bf 0} &0& $n$ & $n$ & $n$& 0 & 2 & final, card=14 \\ \hline
$\rho_{212}=$ & 1 &{\bf 1} &0 &0 & 2 & 2& 0 & 2 & final, card=8 \\ \hline
$\rho_{213}=$ & 1 & 0 &{\bf 1}& 0 & $n$ & $n$ & 0 & 2 & final, card=6 \\ \hline

 \end{tabular}

 \vspace{3mm}
{\sl Table 5: Sketching the trimmed $n$-algorithm}

\vspace{5mm}
\noindent And so it goes on. In the end we find that exactly  the 110 sets 
$$X\in \rho_{111}\uplus \rho_{112}\uplus\rho_{122}\uplus\rho_{1211}\uplus\rho_{1212}\uplus\rho_{22}\uplus\rho_{211}\uplus\rho_{212}
\uplus\rho_{213}$$
yield triangle-free graphs $G_6[X]$. Since $G_6$ is chordal, here "triangle-free" amounts to "forest".

\begin{center}
\includegraphics[scale=1.1]{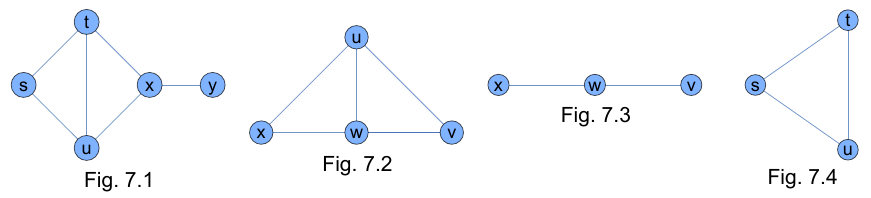}
\end{center}

\vspace{3mm}
{\bf 10.2} In Section 6 the constraints $Y\in{\cal S}$ were imposed one by one. Suppose now the noncover $n$-algorithm must be applied to a set family ${\cal S}\s{\cal P}(W)$ where $|Y|=2$ for all $Y\in{\cal S}$. Then the noncovers wrt $\cal S$ can be viewed as the anticliques  of the graph $G'=(W,{\cal S})$; so $Mod({\cal S})=Acl(G')$. Moreover, we can suitably bundle constraints and then impose whole bundles at once. To fix ideas, consider the graph $G'$ in Figure 4B (ignore its label $Aux(G_4)$) and the bundle ${\cal B}_6$ of all edges incident with $6$, thus
$${\cal B}_6:=\big\{\{6,1\},\{6,4\},\{6,7\},\{6,8\},\{6,10\}\big\}\s {\cal S}.$$
Then $Y\s W=\langle 16\rangle$ is a noncover wrt ${\cal B}_6$ iff $6\not\in Y\ or\ (6\in Y\ and\ 1,4,7,8,10\not\in Y)$.
This gives rise to the (a,c)-wildcard (as in Table 4 blanks and  2's are the same thing):

\vspace{3mm}
\begin{tabular}{c|c|c|c|c|c|c|c|c|c|c|c|c|c|c|c|c|c|l} 
 & 1 & 2 & 3 & 4 & 5 & 6 & 7 & 8 & 9 &
 10& 11& 12&13 & 14&15&16&17& \\ \hline
 &   &    &  &   &   &   &   &   & & &   &  & &   &   &   & & \\ \hline
$r:=$  & $c$ &  &  & $c$ &  & $a$ & $c$ & $c$ &  & & & & & & & && \\ \hline
$r_1:=$  & $\bf 2$ &  &  & $\bf 2$ &  & $\bf 0$ & $\bf 2$ & $\bf 2$ &  & & &  &  & & & &   & \\ \hline
$r_2:=$  & $\bf 0$ &  &  & $\bf 0$ &  & $\bf 1$ & $\bf 0$ & $\bf 0$ &  & & & & & & & &  & \\ \hline
 &   &    &  &   &   &   &   &   & & &   &  & &   &   &   & & \\ \hline
 $\ol{r}:=$  & 2 & 0 & 0 & 1 & $a$ & 0 & 0 & 0 & $c$ & 1& 2& 1&2 & 0& 0&0 &2& \\ \hline
  
   $r':=$  & 2 & 0 & 0 & 1 & $\bf 1$ & 0 & 0 & 0 & $\bf 0$ & 1& 2& 1&2 & 0& 0&0 &2& \\ \hline
    $r^*:=$  & 2 & 0 & 0 & 1 & 1 & 0 & 0 & 0 & 0 & 1& $\bf 1$& 1&$\bf 1$ & 0& 0&0 &2& \\ \hline
\end{tabular}

\vspace{2mm}
{\sl Table 6: Glimpsing the (a,c)-algorithm}

\vspace{2mm} The definition of the (a,c)-wildcard achieves that $r=Mod({\cal B}_6)$.  By imposing\footnote{This is achieved by the (a,c)-algorithm glimpsed in 7.4.} all bundles ${\cal B}_1,...,{\cal B}_{17}$ yields $Mod({\cal S})=Acl(G')$. One of the final (a,c)-rows is $r^*$ in Table 6.

\vspace{2mm}
{\bf 10.2.1} Let us disclose why the proper name for the graph from Fig. 4B is $Aux(G_4)$. Recall from the proof of Theorem 9 that imposing all type 2 constraints $[\{\alpha,\beta\}]$ amounts to applying the noncover $n$-algorithm to a family of 2-elements sets. One checks that the latter match the edges of $Aux(G_4)$. For instance $[\{6,7\}]$ is a type 2 constraint wrt $G_4$ because the path determined by the edges $6$ and $7$ (i.e. the path $(b,f,h)$) is chordless. And indeed, 
$[\{6,7\}]$ being a type 2 constraint matches the edge between the vertices $6$ and $7$ of $Aux(G_4)$. As another illustration, in $Aux(G_4)$ neither $6$ and $11$, nor $6$ and $12$, are adjacent. This is because neither pair yields a type 2 constraint.
(The fact that $6,12$ {\it do} occur in the type 1 constraint $[\{6,12,9\}]$ is irrelevant.)

\vspace{2mm} Therefore the sets collected in the final 012ac-rows (one of which being $\ol{r}$ in Table 6) satisfy all type 2 constraints. It remains to impose all type 1 constraints. As we know, they match the triangles of $G_4$, which we collect here
$$(24)\quad I=\{1,9,15\},\ II=\{2,5,15\},\ III=\{3,11,14\},\ IV=\{4,5,12\},\ V=\{4,10,11\}$$
$$VI=\{5,11,13\},\ VII=\{6,9,12\},
\ VIII=\{7,8,17\},\ IX=\{8,10,16\},\ X=\{10,12,13\}$$
Let us impose the constraints I to X in (24) upon $\ol{r}$. The constraints $II,III,VIII,IX$ happen to hold for all sets $Y\in\ol{r}$ already. For instance $II$ holds because $2,15\not\in Y$ for all $Y\in\ol{r}$. Each of the remaining constraints fails in at least one
$Y\in\ol{r}$. Note that all $X\in \ol{r}$ satisfy $\{4,12\}\s X$. Therefore, in order to satisfy IV we need to set $a:=1$ in $\ol{r}$. This yields the row $r'$ which happens to satisfy all constraints except V,VI,X. One checks that $r^*$ in Table 6 consists of all bitstrings of $r'$ that also satisfy V,VI,X.

Since all type 2 constraints were secured beforehand, Lemma 8 guarantees that all four members of $r^*$ are Cli-Pacs. The largest Cli-Pac is $Y=\{1,4,5,10,11,12,13,17\}$. In terms of vertices we have $Y=E(\Pi)$ for $\Pi:=\{\{a,b\},\{c,d,f,g\},\{h,i\},\{e\}\}$.

 \vspace{3mm}
{\bf  10.3} Similar to 10.1, where we trimmed the noncover $n$-algorithm without affecting the bound in Theorem 6, here we sketch how the implication $n$-algorithm can be trimmed
without affecting the bound in Theorem 11. 

 We previously used the don't-care "2" and the $n$-wildcard $(n,n,....,n)$ which, recall, means "at least one 0 here". Recall from Section 3 that $(n(2),...,n(2))$ means "at least two 0's here" and that
 $(\gamma, \cdots ,\gamma)$  means "exactly one 0 here''. Thus if $(n(2),...,n(2))$ has length $k$ then $|(n(2),...,n(2))|=2^k-1-k$ because the bitstring $(1,1,...,1)$ and the $k$ bitstrings $(0,1,1...,1),(1,0,1,...,1),(1,1,...,1,0)$ do not belong to 
$(n(2),...,n(2))$.

\vspace{3mm}
Probably the extra wildcards $(n(2),...,n(2))$ and $(\gamma, \cdots ,\gamma)$ pay off the most for sparse graphs with a sizeable amount of chordless cycles, e.g. planar graphs.
For instance $G_7=(V^7,E_7)=(V^7,\{1,2,..,11\})$ in Figure 8  has five chordless cycles of cardinalities 7,7,6,5,4 (check). Upon imposing them ad hoc by hand the author represented ${\cal F}_E(G_7)=Cl(\Sigma^{G_7})$ as disjoint union of 29 multivalued rows, the fattest\footnote{Not so fat, but lovely still: $(\gamma,\gamma,n,n,1,n',n',0,\gamma',\gamma',\gamma')\s {\cal F}_E(G_7)$.} being $\ol{r}:=(n,n,n,n,0,n',n',n',n(2),n(2),n(2))$, which houses $15\cdot 7\cdot 4=420$ members of ${\cal F}_E(G_7)$. For instance $(1,1,1,0,0,0,1,1,0,1,0)\in \ol{r}$ "is" the edge-set
$E_7(\Pi)=\{1,2,3,7,8,10\}$ of the Conn-Pac $\Pi\in {\cal F}_V(G_7)$  whose three parts consist, respectively, of the white, blue and black vertices in Figure 8.

\begin{center}
\includegraphics[scale=0.6]{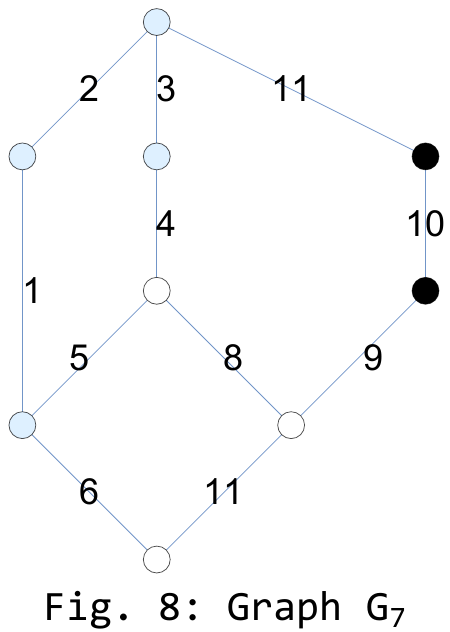}
\end{center}

As to the competitors, the number of 012n-rows and 012-rows required by the standard implication $n$-algorithm and Mathematica's {\tt BooleanConvert} respectively, were 84 and 164. All three methods agreed that $|{\cal F}_E(G_7)|=1190$.

\section{References}

\begin{itemize}
\item[\bf BM] J.A. Bondy, U.S.R. Murty, Graph Theory, Springer 2007.
\item[\bf CH] Y. Crama, P. Hammer (editors),  Boolean functions, Cambridge University Press 2011.
  \item[\bf DCLJ]  E. Dias, D. Castonguay, H. Longo, W. Jradi, Efficient enumeration of chordless cycles,  arXiv:1309.1051v4.
\item[\bf FJ] M. Farber, R.E. Jamison, Convexity in graphs and hypergraphs, SIAM J. Alg. Disc. Math. (1986) 433-444.
\item[\bf M] A. Montina, Output-sensitive algorithm for generating the flats of a matroid,  arXiv:1107.4301v1.

\item[\bf S] A. Schrijver, Combinatorial Optimization, Springer 2003.
\item[\bf SA] A.R. Sharafat, O.R. Arouzi, RECURSIVE CONTRACTION ALGORITHM: A NOVEL AND EFFICIENT GRAPH TRAVERSAL METHOD FOR SCANNING ALL MINIMAL CUT SETS, Iranian Journal of Science and Technology, Transaction B, Engineering, Vol. 30 (2006) 749-761.
    \item[\bf W1] M. Wild, Compactly generating all satisfying truth assignments of a Horn formula, Journal on Satisfiability, Boolean Modeling and Computation 8 (2012) 63-82.
    \item[\bf W2] M. Wild, Compression with wildcards: All metric induced subgraphs, arXiv.2409.08363v3.
     \item[\bf W3] M. Wild, Enumerating all geodesics, submitted.
      \item[\bf W4] M. Wild, Compression with wildcards: All anticliques (or all maximum ones), submitted.
    
\end{itemize}

\end{document}